\documentclass[reqno]{amsart}

\newif\ifmydraft
\mydraftfalse

\usepackage{a4wide}
\usepackage{subfigure}
\usepackage{amssymb}
\usepackage{color}
\usepackage{url}
\usepackage{mathtools}
\usepackage{ifpdf}
\usepackage{enumerate}
\usepackage[font=small,labelfont=bf]{caption}
\usepackage{framed}

\ifpdf
   \usepackage[pdftex]{graphicx}
   \usepackage[bookmarksopen=false,%
               bookmarks,%
               pdftex,%
               colorlinks]{hyperref}
\else
  
  \usepackage{graphicx}
  \ifmydraft
    \usepackage[colorlinks]{hyperref}
    \usepackage[color]{showkeys} 
  \fi
\fi

\DeclareMathOperator{\Ker}{Ker}

\DeclareMathOperator{\ad}{ad}
\DeclareMathOperator{\characteristic}{char}
\DeclareMathOperator{\length}{length}
\newcommand{\Aa}{{A_{n - 1}}}
\newcommand{\Bb}{{B_{n - 1}}}
\newcommand{\Cc}{{C_{n/2}}}
\newcommand{\Dd}{{D_n}}
\newcommand{\FA}{\Ff{\GA}}
\newcommand{\FB}{\Ff{\GB}}
\newcommand{\FC}{\Ff{\GC}}
\newcommand{\FD}{\Ff{\GD}}
\newcommand{\FFD}{\mathbb{F}_{D;n}}
\newcommand{\FF}{\mathbb{F}}
\newcommand{\Ff}[1]{{\mathcal{F}_{#1}}}
\newcommand{\Fff}{\mathcal{F}}
\newcommand{\Fg}{\Ff{\Gamma}}
\newcommand{\GA}{{\Gamma_{A;n}}}
\newcommand{\GB}{{\Gamma_{B;n}}}
\newcommand{\GC}{{\Gamma_{C;n}}}
\newcommand{\GD}{{\Gamma_{D;n}}}
\newcommand{\LL}{\mathcal{L}}
\newcommand{\Ll}[2]{{\LL (#1, #2)}}
\newcommand{\MA}{\MM_\GA}
\newcommand{\MB}{\MM_\GB}
\newcommand{\MC}{\MM_\GC}
\newcommand{\MD}{\MM_\GD}
\newcommand{\MM}{\mathcal{M}}
\newcommand{\None}{\mathbb{N}_+}
\newcommand{\PPeeennnn}[4]{f (#1, #2)#1 #3 #4}
\newcommand{\PPeeennn}[4]{f (#1, #3)#1 #2 #4}
\newcommand{\PPeeenn}[4]{f (#1, #2 #3)#1 #4}
\newcommand{\Peeennnn}[4]{f (#1, #2)#4 #1 #3}
\newcommand{\Peeennn}[4]{f (#1, #3)#4 #1 #2}
\newcommand{\Peeenn}[4]{f (#1, #2 #3) #4 #1}
\newcommand{\Peennnn}[3]{f (#1, #2) #1 #3}
\newcommand{\Peennn}[3]{f (#1, #3) #1 #2}
\newcommand{\Peenn}[3]{f (#1, #2 #3) #1}
\newcommand{\SL}{\mathrm{SL}}
\newcommand{\Sp}{\mathrm{Sp}}
\newcommand{\VV}{\mathcal{V}}
\newcommand{\Vv}[1]{{\VV (#1)}}
\newcommand{\XA}{\Xx{\GA}}
\newcommand{\XB}{\Xx{\GB}}
\newcommand{\XC}{\Xx{\GC}}
\newcommand{\XD}{\Xx{\GD}}
\newcommand{\Xx}[1]{{X (#1)}}
\newcommand{\abs}[1]{\lvert#1\rvert}
\newcommand{\card}[1]{\left\lvert#1\right\rvert}
\newcommand{\downup}{\mathord{\downarrow\!\Uparrow}}
\newcommand{\down}{\mathord{\downarrow}}
\newcommand{\ff}{\mathfrak{f}}
\newcommand{\gll}{\mathfrak{gl}}
\newcommand{\hh}{\mathfrak{h}}
\newcommand{\inref}[1]{{}\stackrel{\eqref{#1}}{\in}{}}
\newcommand{\isref}[1]{{}\stackrel{\eqref{#1}}{=}{}}

\newcommand{\monV}[5]{x_{#1} x_{#2} x_{#3} x_{#4} #5}
\newcommand{\oo}{\mathfrak{o}}
\newcommand{\psiA}{\psi_\GA}
\newcommand{\psiB}{\psi_\GB}
\newcommand{\psiC}{\psi_\GC}
\newcommand{\psiD}{\psi_\GD}
\newcommand{\psiG}{\psi_\Gamma}
\newcommand{\ranF}{\rangle_\FF}
\newcommand{\ranGp}{\rangle_\mathrm{Gp}}
\newcommand{\ranIdl}{\rangle_\mathrm{Idl}}
\newcommand{\ranLie}{\rangle_\mathrm{Lie}}
\newcommand{\sll}{\mathfrak{sl}}
\newcommand{\spp}{\mathfrak{sp}}
\newcommand{\tensor}{\mspace{-2mu}\otimes\mspace{-2mu}}
\newcommand{\updown}{\mathord{\uparrow\!\Downarrow}}
\newcommand{\up}{\mathord{\uparrow}}

\newcounter{storedequation} 
\let\storedtheequation=\theequation
\newenvironment{lettereqns}[1]{%
  \setcounter{storedequation}{\value{equation}}%
  \setcounter{equation}{0}%
  \renewcommand{\theequation}{#1\arabic{equation}}}{%
  \setcounter{equation}{\value{storedequation}}%
  \renewcommand{\theequation}{\storedtheequation}}

\newcounter{casecounter}
\newenvironment{caselist}{%
  \begin{list}{\textbf{Case \arabic{casecounter}:}}{%
      \usecounter{casecounter}%
      \setlength{\leftmargin}{2mm}%
      \setlength{\listparindent}{0cm}}}{%
  \end{list}}

\theoremstyle{plain}
\newtheorem{theorem}{Theorem}
\newtheorem{lemma}[theorem]{Lemma}
\newtheorem{corollary}[theorem]{Corollary}

\newtheoremstyle{example}%
  {\topsep}
  {\topsep}
  {}
  {}
  {\scshape}
  {.}
  { }
  {}
\theoremstyle{example}
\makeatletter
\newenvironment{myleftbar}{%
  \fram@d {\advance\hsize-\width}}%
 {\endfram@d}
\makeatother
\newcommand{\newexample}[2]{
  \newtheorem{x#1}[theorem]{#2}
  \newenvironment{#1}{%
    \begin{myleftbar}
      \begin{x#1}
      }{%
      \end{x#1}
    \end{myleftbar}
  }
}
\newexample{notation}{Notation}
\newexample{definition}{Definition}
\newexample{example}{Example}

\numberwithin{equation}{section}
\numberwithin{theorem}{section}
\numberwithin{figure}{section}
\numberwithin{table}{section}

\makeatletter
\renewcommand{\p@subfigure}{\thefigure}
\makeatother

\author[J.C.H.W.~in~'t~panhuis]{Jos~in~'t~panhuis}
\address[Jos~in~'t~panhuis]{
Department of Mathematics and Computer Science\\
Technische Universiteit Eindhoven\\
P.O. Box 513, 5600 MB Eindhoven, Netherlands}
\email{j.c.h.w.panhuis@tue.nl}

\author[E.J.~Postma]{Erik~Postma}
\address[Erik~Postma]{
Department of Mathematics and Computer Science\\ 
Technische Universiteit Eindhoven\\ 
P.O. Box 513, 5600 MB Eindhoven, Netherlands} \email{e.j.postma@gmail.com}

\author[D.A.~Roozemond]{Dan~Roozemond}
\address[Dan~Roozemond]{
Department of Mathematics and Computer Science\\
Technische Universiteit Eindhoven\\
P.O. Box 513, 5600 MB Eindhoven, Netherlands}
\email{d.a.roozemond@tue.nl}

\title{Extremal Presentations for Classical Lie Algebras}

\begin{document}

\maketitle

\begin{abstract} The long-root elements in Lie algebras of Chevalley
type have been well studied and can be characterized as extremal
elements, that is, elements $x$ such that the image of $(\ad x)^2$
lies in the subspace spanned by $x$. In this paper, assuming an
algebraically closed base field of characteristic not $2$, we find
presentations of the Lie algebras of classical Chevalley type by means
of minimal sets of extremal generators. The relations are described by
simple graphs on the sets. For example, for $C_n$ the graph is a path
of length $2n$, and for $A_n$ the graph is the triangle connected to a
path of length $n-3$. \end{abstract}

\section{Introduction}

\begin{figure}
  \centering%
  \begin{minipage}{0.45\linewidth}
    \centering%
   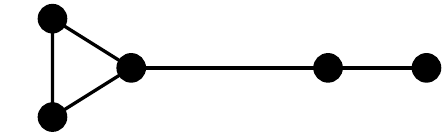
    \caption{The graph $\GA$ for $\sll_n$.}
    \label{fig:twt}%
  \end{minipage} %
  \begin{minipage}{0.45\linewidth}
    \centering%
    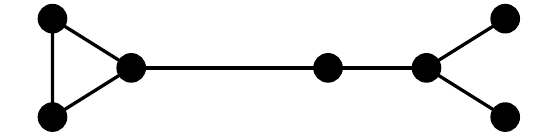
    \caption{The graph $\GB$ for $\oo_{2n-1}$.}
    \label{fig:trspt}%
  \end{minipage}
  \\
  \begin{minipage}{0.45\linewidth}
    \centering%
    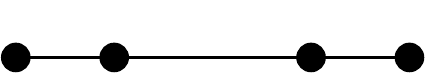
    \caption{The graph $\GC$ for $\spp_n$.}%
    \label{fig:linear}%
  \end{minipage} %
  \begin{minipage}{0.45\linewidth}
    \centering%
    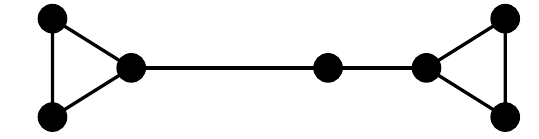
    \caption{The graph $\GD$ for $\oo_{2n}$.}%
    \label{fig:twotr}%
  \end{minipage} 
\end{figure}

A nonzero element $x$ of a Lie algebra $\LL$ over a field $\FF$ of
characteristic not $2$ is called \emph{extremal}\index{extremal
element} if $[x, [x, \LL]] \subseteq \FF x$. Extremal elements are a
well-studied class of elements in simple finite-dimensional Lie
algebras of Chevalley type: they are the long root elements.
In~\cite{cit:cosuw01}, Cohen, Steinbach, Ushirobira and Wales have
studied Lie algebras generated by extremal elements, in particular
those of Chevalley type. The authors also find the minimum size of a
set of generating extremal elements for the Lie algebras of Chevalley
type and find such minimal generating sets of extremal elements
explicitly. In the present paper, we also find such minimal generating
sets of extremal elements explicitly for the four classical families
of Lie algebras: those of type $A_n$, $B_n$, $C_n$ and $D_n$. We will
do this in a more geometrical setting and will find criteria for sets
of extremal elements to generate Lie algebras of this type.

By Lemma~\ref{lm:2extelts}, each Lie algebra generated by a pair of
linearly independent extremal elements is in one of only three
isomorphism classes: either the two-dimensional commutative Lie
algebra, or the so-called Heisenberg Lie algebra $\hh$, or $\sll_2$.
Given a generating set $S$ of extremal elements, we examine the
subalgebras generated by pairs of these elements. These give rise to
graphs: the vertices correspond to the elements of $S$, and two
vertices are adjacent if the corresponding extremal elements generate
a three-dimensional algebra and nonadjacent if they commute.  We will
say that the Lie algebra generated by $S$ \emph{realizes} this graph.

Following experiments using the GAP computer algebra system 
\cite{cit:GAP} and the GBNP package \cite{cit:GBNP} we conjectured one such graph
for \emph{each} Lie algebra of classical Chevalley
type, depicted in Figures~\ref{fig:twt} up to~\ref{fig:twotr}.

In this paper we 
show that if a Lie algebra realizes one of the graphs in 
Figures~\ref{fig:twt} up to~\ref{fig:twotr}, then in the
generic case it is isomorphic to the Lie algebra of the corresponding
Chevalley type, in the following sense. Given a graph $\Gamma$, we
define a vector space $\Vv{\Gamma}$ parametrizing the Lie algebras
that realize $\Gamma$. Let $\Xx{\Gamma}$ be the subset of
$\Vv{\Gamma}$ of values $\ff$ for which the associated Lie algebra
$\Ll{\Gamma}{\ff}$ has maximal dimension among such algebras for the
same graph $\Gamma$. We will see in Lemma~\ref{thm:x-is-a-variety}
that $\Xx{\Gamma}$ carries the structure of an affine variety. The
following theorem and analogues for the three other families of
Chevalley type Lie algebras will be our main results.
\begin{theorem}
  Let $n \geq 5$. Let $\GD$ be the graph of Figure~\ref{fig:twotr}.
  There is an open dense subset $S$ of $\Xx{\GD}$ such that, if $\ff
  \in S$, then $\Ll{\GD}{\ff}$ is isomorphic to the Lie algebra of
  type $\Dd$.
\end{theorem}

The graphs~$\GA$, $\GB$, and $\GC$ are subgraphs of $\GD$, and as such
induce extra commuting relations. Therefore, the $\Dd$ case is the most
complicated one, from which we will deduce the conclusions for the
other cases.

\subsection{Contents and Strategy}
In the rest of this section we will introduce some conventions and
notation that we will use in this paper. In
Section~\ref{sec:preliminaries} we review some of the underlying
theory. In Section~\ref{sec:general-framework} we show how to work
with abstract Lie algebras that realize a given graph. We apply this
to four proposed graphs~$\GA$, $\GB$, $\GC$ and~$\GD$ for the classical Chevalley types in
Section~\ref{sec:monomials}.
In Section~\ref{s:parameter-space}, we
extend this study to the parameter space $\VV$ referred to above.
In Theorems~\ref{thm:realizeslinear},
\ref{thm:realizestwt}, \ref{thm:realizestwotr} and~\ref{thm:realizestrspt},
we give concrete realizations of the Lie
algebras of types~$A_n$, $B_n$, $C_n$ and~$D_n$ corresponding to the
graphs from Figures~\ref{fig:twt} up to~\ref{fig:twotr}.

Finally, in Section~\ref{s:alwaysours} we prove the main results of
this paper: we show that a Lie algebra $\LL$ realizing one of the
graphs~$\GA$, $\GB$, $\GC$ and~$\GD$ is in the generic case a
quotient of the realization $\MM$ found in Section~\ref{s:realizations}.
Since $\MM$ is simple in most cases, it will follow that $\LL$ and
$\MM$ are isomorphic. The only exception is $A_n$ if $p \mid n+1$.

This paper was inspired by the Masters thesis of the third
author~\cite{cit:rooze05} and reported on more extensively in
the second author's Ph.D.~thesis~\cite{cit:postm07}.

\subsection{Conventions and notation}
For the rest of this paper, $\FF$ will be an algebraically closed
field of characteristic not $2$ and $\LL$ will be a Lie algebra over
$\FF$.

Since we approach the matter from the angle of the generating sets of
abstract extremal elements, we let $n$ be the number of generating
extremal elements. In Theorems
Section~\ref{s:realizations}, this will mean
that we study, for example, the Lie algebra of type $\Cc$, defined
over a vector space of dimension $n$, where we have to assume that $n$
is even. In that section, it would be more convenient to study the Lie
algebra of type $C_n$ instead, but we choose consistency over
convenience and keep the meaning of $n$ as the number of extremal
generators.

\medskip

  If no confusion is possible, we write $xy$ for $[x, y]$, and $xyz$
  for $[x, [y, z]]$; we will write $(xy)z$ for $[[x, y], z]$. So,
  anticommutativity and the Jacobi identity will be written as
  \begin{align}
    xx &= 0 \label{e| L} \tag{AC} 
  \end{align}
and
  \begin{align}
		xyz + yzx + zxy &= 0. \label{e| J} \tag{J}
	\end{align}

\medskip

  We will often work with long products of indexed elements. We use
  the following notation to make these products somewhat manageable.
  The general idea is that we put two numbers in the subscript with an
  operator consisting of one or two arrows in between, such as $x_{5
    \updown 2}$; the first factor in the product is then indexed by
  the first number, after which we iterate adding (for up arrows) or
  subtracting (for down arrows) one (for single stroke arrows) or two
  (for double stroke arrows) to the index until we encounter the last
  number, where every step gives the next factor for this product. So
  the previous example $x_{5 \updown 2}$ is short for $x_5 x_6 x_4 x_5
  x_3 x_4 x_2$.

  In particular, there are four operators that we use, defined more
  precisely as follows.
   If $i \leq j$, the notation $x_{i \up j}$
  will mean $x_i x_{i + 1} x_{i + 2} \dotsm x_{j - 1} x_j$, and $x_{j
    \down i}$ will mean $x_j x_{j - 1} x_{j - 2} \dotsm x_{i + 1}
  x_i$.  Furthermore, $x_{j \updown i}$ will mean $x_j x_{j + 1} x_{j
    - 1} x_j x_{j - 2} x_{j - 1} \dotsm x_{i + 1} x_{i + 2} x_i$, and
  similarly, $x_{i \downup j}$ will mean $x_i x_{i - 1} x_{i + 1} x_i
  x_{i+2} x_{i+1} x_{i + 3} \dotsm x_{j - 1} x_{j - 2} x_j$.

  We will also use constructions such as $x_{3 \up 6} x_{4 \updown
    1}$, which will mean $x_3 x_4 x_5 x_6 x_4 x_5 x_3
  x_4 x_2 x_3 x_1$.  Occasionally, it will be convenient to include in
  a set of monomials of the form, say, $x_{j \down i} x_{i - 2}$ the
  case $j = i - 1$; this monomial will then simply be $x_{i - 2}$. So
  in this case $x_{j \down i}$ cannot be seen as a separate monomial.

  We extend the notation to cover the case where we have a sequence
  $i_1, \dotsc, i_k$ of indices: then we write $x_{i_{k \down 1}}$ for
  $x_{i_k} x_{i_{k - 1}} \dotsm x_{i_2} x_{i_1}$.

\medskip

  We say that a set of Lie algebra elements $\{x_i \mid i \in V\}$
  \emph{realizes}\index{realizing a graph} a given graph $\Gamma = (V,
  E)$ if:
  \begin{itemize}
  \item each $x_i$ is an extremal element of $\langle x_j \mid j \in I
    \ranLie$;
  \item vertices $i$ and $j$ are connected if and only if $x_i$ and
    $x_j$ do not commute.
  \end{itemize}
  We will sometimes also say that the Lie algebra $\langle x_i
  \ranLie$ realizes $\Gamma$.
Later in this paper it will be essential that each $x_i$ is nonzero,
which is implied by it being extremal.

\subsection{Related results}
In \cite{cit:drpa08} Lie algebras realized by simply laced affine
Dynkin diagrams were considered. There it was shown that in the
generic case the Lie algebra is of the corresponding finite type. The
diagram for $A_n$ given there can be transformed into Figure
\ref{fig:twt} using a procedure similar to that described in Lemma
\ref{lm:fixtriangle}. The diagram for $D_n$ given in \cite{cit:drpa08}
is related in a less straightforward manner, since Figure
\ref{fig:twotr} has $n$ vertices whereas its affine Dynkin diagram has
$n+1$ vertices.

Although the generators arising from our graphs and some of the
Chevalley generators \cite{cit:cart72} are similar in the sense that
they both correspond to long root elements, no direct relation is
apparent.

\section{Preliminaries}\label{sec:preliminaries}
In this section, we will introduce a bilinear form defined on all Lie
algebras generated by extremal elements, and recall some of its
properties. None of these results are new; most can be found in
e.g.~\cite{cit:cosuw01} and thus we will omit most of the proofs.
We will start by introducing a related family of linear functionals.

For extremal $x$, let $f_x \colon \LL \to \FF$ be the linear map
defined by $xxy = f_x (y) x$. Since $[\cdot, \cdot]$ is bilinear, this
is indeed a linear map. We call $f_x$ the \emph{extremal
  functional}\index{extremal functional} on $x$.

\begin{lemma} \label{lm:f-symm}
  $f_x (y) = f_y (x)$ for all extremal $x,y\in\LL$.
\end{lemma}
\begin{lemma} \label{lm:2extelts}
  Let $\LL = \langle x, y \ranLie$ with $x$ and $y$ extremal and
  linearly independent. Then $\LL$ is isomorphic to the
  two-dimensional commutative Lie algebra, the Heisenberg algebra, or
  $\sll_2$.
\end{lemma}
\begin{lemma} \label{lm:EElieisEEspan}
  If $\LL$ is generated by extremal elements, then it is linearly
  spanned by extremal elements.
\end{lemma}
\begin{lemma}
  If $\LL$ is generated by extremal elements, the definition of $f_x
  (y)$ can be extended to a unique bilinear form $f (x, y)$ on $\LL$
  with $f (x, y) = f_x (y)$ if $x$ is an extremal element.
  This bilinear form is associative and symmetric:
  \begin{equation}
    \forall x,y,z \colon f (x, yz) = f (xy, z)\label{e| f} \tag{AS}
  \end{equation}
  \begin{equation}
  	 \forall x, y \colon f (x, y) = f (y, x). \label{e| ff} \tag{SM}
  \end{equation} 
\end{lemma}

\begin{lettereqns}{P}
  We call $f$ the \emph{extremal form}\index{extremal form}. 
  We will use the following identities involving the extremal form,
  the first two of which go back to Premet and were first used
  in~\cite{cit:chern89}:
  \begin{lemma}\label{lem:premet}
    If $x,y,z\in\LL$ and $x$ extremal, then
	\begin{equation}
      2 (x y)x z={} f(x, y z) x+f(x, z) x y -f(x, y) x z,\label{e| P1}
	\end{equation}
	\begin{equation}
      2 x y x z={} f(x, y z) x-f(x, z) x y -f(x, y) x z,\label{e| P2}
	\end{equation}
	\begin{equation}
      f (x, y x z)={} - f (x, z) f (x, y). \label{e| P5}		
	\end{equation}
  \end{lemma}
  \begin{proof} By the Jacobi identity,
    \begin{align*}
      (xy)xz \isref{e| J} & ((xy)x) z + x (xy) z = - f (x, y) xz + x
      (xy) z, \\
      \intertext{and similarly,}
      (xy)xz \isref{e| L} & -(xz)xy \isref{e| J} - ((xz)x)y - x (xz) y
      \isref{e| J} f (x, z) xy - x x z y - x (x y) z.
    \end{align*}
    Adding these two equations and applying anti-commutativity a few
    times, we obtain Eq.~\eqref{e| P1}. Then we find Eq.~\eqref{e| P2}
    as follows:
    \[
    2 x y x z \isref{e| J} 2 (x y) x z + 2 y x x z \isref{e| P1} f(x, y
    z) x+f(x, z) x y -f(x, y) x z - 2 f (x, z) x y.
    \]
    For Eq.~\eqref{e| P5}, we need the next lemma. The equation
    is then easily obtained as follows:
    \[
    f (x, y x z) \isref{e| f} f (xy, xz) \isref{e| L} - f (yx, xz)
    \isref{e| f} - f (y, xxz) = - f (x, y) f (x, z). \tag*{\qedhere}
    \]
  \end{proof}
\end{lettereqns}

\section{The general framework}\label{sec:general-framework}
In this section, we will establish a framework for dealing with Lie
algebras generated by extremal elements where we prescribe the values
of the extremal form. In the end, we prove
Theorem~\ref{thm:x-is-a-variety} which shows that a certain parameter
space for these prescribed values is an algebraic variety. The main
objective of this section is to introduce the techniques for proving
that theorem. In Section~\ref{s:parameter-space} we will use those
techniques to prove similar theorems for a substantially smaller
parameter space, but then for specific Lie algebra families.

Let $n \in \None$ be fixed.  Let $\Gamma$ be a graph on $n$ numbered
vertices. Let $\Fff$ be the free Lie algebra over $\FF$ on $n$
generators $x_1, \dotsc, x_n$ with the standard grading. We will
construct a quotient of $\Fff$ where the projections of $x_i$ in the
quotient are extremal generators. Let
\[
\Fg = \Fff / \langle x_i x_j \mid \{i, j\} \not\in E(\Gamma)\ranIdl.
\]
$\Fg$ inherits the grading of $\Fff$; this is possible because the
ideal that is divided out is \emph{homogeneous}\index{homogeneous!with respect to a grading} with respect to the grading of $\Fff$, in the
sense that it is spanned by its intersections with the homogeneous
components of $\Fff$.

Let $\ff = (\ff_1, \dotsc, \ff_n)$ be an element of $\Vv{\Gamma}
\coloneqq (\Fg^*)^n$, so it consists of $n$ functionals in the dual of
$\Fg$; we will make sure that $\ff_i$ is the extremal functional
$f_{x_i}$ in the Lie algebra we will construct. To that end, define
the ideal
\[
I_{\Gamma, \ff} = \langle x_i x_i y - \ff_i (y) x_i \mid
y \in \Fg, 1 \leq i \leq n \ranIdl.
\]
When taking $\ff = 0$, we see that $I_{\Gamma, 0}$ is homogeneous with
respect to the standard grading of $\Fg$. Let $\Ll{\Gamma}{\ff} = \Fg
/ I_{\Gamma, \ff}$ and let $\xi_\ff \colon \Fg \to \Ll{\Gamma}{\ff}$
be the natural projection.  We will sometimes omit $\xi_\ff$ if that
does not stand in the way of clarity.  Clearly each $x_i$ is either an
extremal element of $\Ll{\Gamma}{\ff}$ or zero, and $\ff_i (y) = f
(x_i, y)$.  We find the following slight extension of Lemma~4.3
of~\cite{cit:cosuw01}:
\begin{lemma} \label{lm:csuw-extended}
  There is a finite list $\MM_\Gamma$ of monomials in $x_1, \dotsc, x_n$
  satisfying the following properties:
  \begin{enumerate}
  \item $\xi_0 (\MM_\Gamma)$ is a basis of $\Ll{\Gamma}{0}$,
    \label{item:mm_gamma-basis}
  \item if $x_i m \in \MM_\Gamma$, then $m \in \MM_\Gamma$,
    \label{item:xi-m-m}
  \item $\MM_\Gamma$ contains all generators $x_i$,
    \label{item:mm_g-all-xi} and
  \item $\Ll{\Gamma}{\ff} = \langle \xi_\ff (\MM_\Gamma)\ranF$ for
    all $\ff \in \Vv{\Gamma}$. \label{item:spans-all-f}
  \end{enumerate}
\end{lemma}
\begin{proof}
  $\Ll{\Gamma}{0}$ is a quotient of $\LL_n \coloneqq \Ll{K_n}{0}$,
  where $K_n$ is the complete graph. By Theorem~1 of
  Zel'manov~\cite{cit:zel80} (or for characteristic $3$, by Theorem~1
  of Zel'manov and Kostrikin~\cite{cit:zelko90}), we know that $\LL_n$
  is finite-dimensional.  Hence we can find a finite set $\MM_\Gamma$
  of monomials in $\Fg$ satisfying
  conditions~\ref{item:mm_gamma-basis}, \ref{item:xi-m-m}
  and~\ref{item:mm_g-all-xi}, by the following procedure. We start by
  setting $\MM_\Gamma$ equal to $\{x_1, \dotsc, x_n\}$. This set is
  linearly independent because the free Abelian Lie algebra on $n$
  generators is a quotient of $\Ll{\Gamma}{0}$, and the images of the
  $x_i$ in it are linearly independent.  Then we perform a number of
  rounds as follows. In each round we form the monomials $x_i m$,
  where $x_i$ iterates over the generators of $\Fg$ and $m$ iterates
  over the longest monomials in $\MM_\Gamma$ so far. We select a
  subset of these such that its images under $\xi_0$ in $\Ll{\Gamma}{
    0}$ are linearly independent of each other and of the images under
  $\xi_0$ of the elements in $\MM_\Gamma$ so far, and add it to
  $\MM_\Gamma$. Then we continue with the next round if we have added
  any new monomials this round. Since $\Ll{\Gamma}{0}$ is
  finite-dimensional, this procedure terminates after finitely many steps.

  Let $U = \langle \MM_\Gamma \ranF \subset \Fg$.  We now prove
  condition~\ref{item:spans-all-f} by showing that $\xi_\ff (U) =
  \Ll{\Gamma}{\ff}$ for all $\ff \in \Vv{\Gamma}$.  Note that
  $I_{\Gamma, 0}$ is spanned by elements of the form $x_{i_{k \down 1}}
  x_{i_1} r$, with $r$ a monomial in $\Fg$.

  Clearly $\xi_\ff (U) \subset \Ll{\Gamma}{\ff}$. Suppose that it is
  a proper subset; then there are monomials $s \in \Fg$ such that
  $\xi_\ff (s) \notin \xi_\ff (U)$, whence $s \notin U$. Let $s$ be
  such a monomial of lowest degree. Since $\Fg = U + I_{\Gamma, 0}$,
  it is possible to express $s$ as a linear combination of monomials
  in $\MM_\Gamma$ and monomials of the form $x_{i_{k \down 1}} x_{i_1}
  r$. All these monomials have the same degree, because only the
  Jacobi identity and anticommutativity can be used for rewriting, in
  addition to homogeneous elements being $0$. Let $t$ be a monomial of
  the form $x_{i_{k \down 1}} x_{i_1} r$ such that $\xi_\ff (t)
  \not\in \xi_\ff (U)$ and let $t_0 = x_{i_{k \down 1}}$.  Then
  \[
  \ff_{i_1} (r) \xi_\ff (t_0) = \xi_\ff (t) \notin \xi_\ff (U),
  \]
  so $\ff_{i_1} (r) \not= 0$ and $t_0 \notin U$. Since $\deg t_0 <
  \deg t = \deg s$, we have a contradiction. Hence $\xi_\ff (U) = \Ll{\Gamma}{
    \ff}$.
\end{proof}
Define $U = \langle \MM_\Gamma \ranF$ as in the preceding proof. Note
that $\Fg = U + I_{\Gamma, \ff}$ for all $\ff$, not just for $\ff =
0$. Define $I$ and $m_i$ by letting $\MM_\Gamma = \{m_i \mid i \in
I\}$.
\begin{lemma} \label{lm:n_m-exists}
  For every monomial $m \in \Fg$, there exists a map $n_m \colon
  \Vv{\Gamma} \to U$, such that $n_m (\ff) = m \pmod{I_{\Gamma, \ff}}$
  for all $\ff \in \Vv{\Gamma}$ and the following property holds.  If
  $n_m (\ff) = \sum_{i \in I} \alpha_{m, i, \ff} m_i$, then
  $\alpha_{m, i, \ff}$, when regarded as a function in $i$ and $\ff$,
  is a polynomial function in the values of $\ff$ at monomials of
  degree less than $\deg m$.
\end{lemma}
\begin{proof}
  Let $m = x_{i_{k \down 1}}$ be a monomial in $\Fg$ of degree $k$. If
  $k = 1$, we put $n_m (\ff) = m$. We proceed by induction on $\deg
  m$. Since $\Fg = U + I_{\Gamma, 0}$, we can write $m$ as the sum of
  an element $u$ of $U$ and an element $w$ of $I_{\Gamma, 0}$; all
  monomials involved have the same degree, because only the Jacobi
  identity and anticommutativity can be used for rewriting, in
  addition to homogeneous elements being $0$. If we prove that $n_w
  (\ff) = w \pmod{I_{\Gamma, \ff}}$ and that its coefficients
  $\alpha_{m,i,\ff}$ satisfy the polynomiality condition, then setting
  \begin{equation}
    n_m (\ff) = u + n_w (\ff)\label{eq:def-nm-ff-gen}
  \end{equation}
  will be sufficient to show that the lemma holds for $m$.

  We may assume that $w$ is a single monomial. So the proof obligation
  reduces to the case where $m \in I_{\Gamma, 0}$. Then there exist
  $r, h \in \None$ such that $i_r = i_{r-1} = h$ and $m$ is thus of
  the form $x_{i_{k \down r+1}} x_h x_h x_{i_{r-2 \down 1}}$. Hence, by
  the induction hypothesis,
  \[
  m = \ff_h (x_{i_{r-2 \down 1}}) x_{i_{k \down r}} = \ff_h
  (x_{i_{r-2 \down 1}}) n_{x_{i_{k \down r}}} (\ff)
  \pmod{I_{\Gamma, \ff}}.
  \]
  We choose
  \begin{equation}
    n_m (\ff) = \ff_h (x_{i_{r-2 \down 1}}) n_{x_{i_{k \down
          r}}} (\ff), \label{eq:def-nm-ff-inI}
  \end{equation}
  so that
  \[
  \alpha_{m, j, \ff} = \ff_h (x_{i_{r-2 \down 1}})
  \alpha_{x_{i_{k \down r}}, j, \ff}.
  \]
  The coefficients $\alpha_{x_{i_{k \down r}}, j, \ff}$ are, by the
  induction hypothesis, polynomials in the values of $\ff$ at
  monomials of degree less than $k - r + 1 < k = \deg m$, so
  equations~\eqref{eq:def-nm-ff-gen} and~\eqref{eq:def-nm-ff-inI}
  define a map satisfying the conditions in the lemma.
\end{proof}
Note that we do not claim that $n_m$ is uniquely determined by these
conditions.  We choose a map $n_*$ as above and extend it to general
elements of $\Fg$ by linearity.

Let $\Xx{\Gamma} = \{\ff \mid \dim \Ll{\Gamma}{\ff} =
\card{\MM_\Gamma}\}$ and let $R\colon \Xx{\Gamma} \to (U^*)^n$ be the map
that restricts a functional to $U$.
\begin{lemma} \label{lm:Risinjective}
  The restriction map $R$ is injective.
\end{lemma}
\begin{proof}
  Let $\ff \in \Xx{\Gamma}$. Then all $m_i \in \MM_\Gamma$ are linearly
  independent in $\Ll{\Gamma}{\ff}$, so $x_i \not\in I_{\Gamma,
    \ff}$.  Let $m$ be a monomial in $\Fg$. We will show that $\ff_i
  (m)$ can be expressed in the values of $\ff_i$ on monomials in
  $M$. If $m \in M$, there is nothing to prove, so assume $m \not\in
  M$. Since
  \[
  x_i x_i m = \ff_i (m) x_i \pmod{ I_{\Gamma, \ff}},
  \]
  and also
  \[
  x_i x_i m = x_i x_i n_m (\ff) = \ff_i (n_m (\ff)) x_i
  \pmod{I_{\Gamma, \ff}},
  \]
  we find that $\ff_i (m) = \ff_i (n_m (\ff))$. Since $n_m (\ff)$ only
  depends on monomials of lower degree than $m$, we see that $\ff_i
  (m)$ can be expressed in the values of $\ff_i$ at monomials of lower
  degree than $m$. By induction on the degree of $m$, it can therefore
  be expressed ultimately in the values of $\ff_i$ on $M$, as we set
  out to prove.

  Let $\ff, \ff' \in \Xx{\Gamma}$ with $R (\ff) = R (\ff')$. Then $\ff_i$
  and $\ff_i'$ agree on $U$, and thus on $\Fg$, for all $i$. Hence
  $\ff = \ff'$.
\end{proof}
\begin{lemma} \label{lm:imRisclosed}
  $R (\Xx{\Gamma})$ is a closed subset of $(U^*)^n$.
\end{lemma}
\begin{proof}
  For all $\ff \in \Vv{\Gamma}$, let the bilinear anticommutative map
  $[\cdot, \cdot]_\ff \colon U \times U \to U$ be determined by
  \[
  [v, w]_\ff = n_{[v, w]} (\ff).
  \]
  If $\ff \in \Xx{\Gamma}$, then
  \begin{enumerate}
  \item $[\cdot, \cdot]_\ff$ is a Lie multiplication (i.e.~it
    satisfies the Jacobi identity),
  \item $[x_i, [x_i, v]_\ff]_\ff = \ff_i (v) x_i$ for all $v \in U$
    and all $i$, \label{item:xixiv}
  \item $[x_i, x_j]_\ff = 0$ if nodes $i$ and $j$ are not connected by
    a line,\label{item:xixj}
  \item the Lie algebra $(U, [\cdot, \cdot]_\ff)$ is generated by
    $x_1, \dotsc, x_n$. \label{item:generated}
  \end{enumerate}
  On the other hand, if all of the above conditions hold for a
  multiplication map $\mu$, then $(U, \mu)$ is a quotient of
  $\Ll{\Gamma}{\ff}$ of the same dimension, and hence isomorphic to
  $\Ll{\Gamma}{\ff}$. But these conditions are all polynomial in the
  values of $\ff$ on $U$: the Jacobi identity and
  conditions~\ref{item:xixiv} and~\ref{item:xixj} are polynomial in a
  straightforward way, and condition~\ref{item:generated} is always
  satisfied: for every $w = x_{i_{k \down 1}} \in \MM_\Gamma$, we have
  $w = [x_{i_k}, [x_{i_{k - 1}}, \dotsb, [x_{i_2}, x_{i_1}]_\ff \dotsb
  ]_\ff]_\ff$, since $n_w (\ff) = w$. So $R (\Xx{\Gamma})$ is given as
  the zero set of a set of polynomial equations; thus it is closed.
\end{proof}
\begin{theorem} \label{thm:x-is-a-variety}
  $\Xx{\Gamma}$ carries a natural structure of an affine variety.
\end{theorem}
\begin{proof}
  The restriction map $R$ is a continuous bijection of $\Xx{\Gamma}$ with
  a Zariski closed subset of $(U^*)^n$.
  Clearly the restriction map $R$ is continuous. The preceding two
  lemmas show that it is injective and that its image is closed.
\end{proof}

\section{The monomials}\label{sec:monomials}
In Figures~\ref{fig:twt} up to~\ref{fig:twotr} we defined four
graphs, $\GA$, $\GB$, $\GC$ and $\GD$, to be used for $\Gamma$ in the
framework of Section~\ref{sec:general-framework}. In this section, we
will construct the basis $\MM_\Gamma$ of $U$ explicitly for each such
$\Gamma$, resulting in Theorems~\ref{t| L1 monomials}, \ref{t| L2
  monomials}, \ref{t| L3 monomials} and~\ref{t| L4 monomials}.

Clearly, each of the algebras $\Fg$ is defined by a subset of the four
following relations.
\begin{lettereqns}{R}
  \begin{align}
    x_i x_j &= 0 & \text{for all $i, j$ with $\abs{i - j} > 1, \{1,
      3\} \not= \{i, j\} \not= \{n - 2, n\}$,} \label{e| R1} \\
    x_1 x_3 &= 0, \label{e| R4} \\
    x_{n - 2} x_n &= 0, \label{e| R2} \\
    x_{n - 1} x_n &= 0. \label{e| R3}
  \end{align}
\end{lettereqns}
We will use the following technical lemma:
\begin{lemma}\label{lm| identities}
  Let $a,b,n\in\None$, $i,j,k,\ell,m,i_1,\dotsc,i_a,j_1,\dotsc,j_b\in
  \{1,\dotsc,n\}$. Let $\Gamma$ be one of the graphs
  $\GA$, $\GB$, $\GC$ and $\GD$, let $\ff \in \Vv{\Gamma}$, 
  let $x_i$ be the standard generators of $\Ll{\Gamma}{\ff}$, and let $t,u \in \Ll{\Gamma}{\ff}$.
  Furthermore, let
  $x_{i_p}$ commute with $x_{j_q}$ for all $p$ and $q$ and let $x_i$
  commute with $x_j$. For Eq.~\eqref{e| Q2} only, assume that $i
  < n-2$.  Then: \pagebreak[1]
  \begin{lettereqns}{Q}
    \begin{align}
      x_j x_{i} t={}&x_i x_{j} t,\label{e| Q1} \\
      x_i x_{i+1} x_{i+2} x_i
      t={}&\frac{1}{2}\left(f (x_i,  x_{i+1} x_{i+2} t)x_i-f (x_i,
        x_{i+2} t)x_i x_{i+1} -f (x_i, x_{i+1}) x_i x_{i+2}
        t\right),\label{e| Q2} \\
      \monV{k}{\ell}{m}{k}{t}={}&\frac{1}{2}\left(\Peeenn{x_{k}}{x_m}{t}{x_{\ell}}
        +\Peeennn{x_{k}}{x_m}{t}{x_{\ell}}-
        \Peeennnn{x_{k}}{x_m}{t}{x_{\ell}}
        +\Peenn{x_{k}}{x_{\ell}}{x_{m} t}\nonumber\right.\\
      &{}-\Peennn{x_{k}}{x_{\ell}}{x_{m} t}-\Peennnn{x_{k}}{x_{\ell}}{x_{m}
        t}-\PPeeenn{x_{k}}{x_{\ell}}{x_m}{t}+
      \PPeeennn{x_{k}}{x_{\ell}}{x_m}{t}\nonumber\\
      &{}+\PPeeennnn{x_{k}}{x_{\ell}}{x_m}{t} -\Peenn{x_{k}}{x_{m}}{x_{\ell}
        t}+\Peennn{x_{k}}{x_{m}}{x_{\ell} t}+
      \Peennnn{x_{k}}{x_{m}}{ x_{\ell} t }\nonumber\\
      &{}+\left.\Peeenn{x_{k}}{x_{\ell}}{t}{x_m}-
        \Peeennn{x_{k}}{x_{\ell}}{t}{x_m}-
        \Peeennnn{x_{k}}{x_{\ell}}{t}{x_m}\right) +
      \monV{k}{m}{\ell}{k}{t},\label{e| Q3} \\ \pagebreak[1]
      f(u,\monV{k}{\ell}{m}{k}{t}))={}&\frac{1}{2}\left(\Peeenn{x_{k}}{x_m}{t}{f (u, x_{\ell}})
        +\Peeennn{x_{k}}{x_m}{t}{f (u, x_{\ell}})-
        \Peeennnn{x_{k}}{x_m}{t}{f (u, x_{\ell}})\nonumber\right.\\
      &{}+f (x_k, x_\ell x_mt)f (u, x_k) -f (x_k, x_mt)f (u, x_kx_\ell)-
      f (x_k, x_\ell)f (u, x_kx_mt) \nonumber\\ &{}- f (x_k, x_\ell x_m)f (u, x_kt)
      +f (x_k, x_m)f (u, x_kx_\ell t)+f (x_k, x_\ell)f (u, x_kx_mt)
      \nonumber\\ &{}- f (x_k, x_mx_\ell t)f (u, x_k)+
      f (x_k, x_\ell t)f (u, x_kx_m)+
      f (x_k, x_m)f (u, x_kx_\ell t) \nonumber\\
      &{}+\left.\Peeenn{x_{k}}{x_{\ell}}{t}{f (u, x_m})-
        \Peeennn{x_{k}}{x_{\ell}}{t}{f (u, x_m})-
        \Peeennnn{x_{k}}{x_{\ell}}{t}{f (u, x_m}) \nonumber\right)\\ &{}+
      f(u,\monV{k}{m}{\ell}{k}{t}) ,\label{e| Q3a} \\ x_{i_1} x_{j_1}
      x_{j_2} \dotsm x_{j_b}={}&0.
      \label{e| Q4} 
    \end{align}
  \end{lettereqns}
\end{lemma}

The proof of this lemma is straightforward, using the identities from Lemma \ref{lem:premet} and the Jacobi identity.

\subsection{Monomials for $\GD$}
In the rest of this section, we will let $\ff \in \Vv{\GD}$ be
given and consider $\Ll{\GD}{\ff}$, that
is, only the relations in~\eqref{e| R1} from page~\pageref{e| R1} are
divided out. In Theorem~\ref{t| L1 monomials}, we will give a list of
$2n^2-n$ monomials in $x_1,\dotsc,x_n$ and prove that they span
$\Ll{\GD}{\ff}$ linearly.

Note that the precise contents of this list is of little general
importance. For example, the number of classes may be reduced using a
more clever notation. However, one needs to fix such a list for the proof of
Theorem \ref{thm:twotralways}, and the one presented below suffices.

\begin{theorem}\label{t| L1 monomials}
  Let $\MD$ be the set consisting of the following
  monomials:
  \begin{align*}
    y^1_{k,m}&{}=x_{k \down m},&n\geq k\geq m\geq 1,\\
    y^2_{k,m}&{}=x_{k \up n-2} x_{n \down m},& n-2\geq k>m\geq 1,\\
    y^3_{k,m}&{}=x_{k \down m+1}x_{m-1 \updown 1}, & n \geq k\geq
    m\geq 3,\\
    y^4_{k,m}&{}=x_{k \up n-2}x_{n \down m + 1} x_{m - 1 \updown 1},&n-2\geq k\geq m\geq 3,\\
    y^5_{m}&{}=x_n x_{n-2 \down m},& n-2\geq
    m\geq 1,\\
    y^6_m&{}=x_{n-1} x_n x_{n-2 \down m},
    & n-2\geq m\geq 1,\\
    y^7_{k}&{}= x_{k \down 3} x_1,&n\geq k\geq 3, \\
    y^8_{k}&{}=x_{k \up n-2}x_{n \down 3} x_1,&n-2\geq k\geq 2,\\
    y^9_{m}&{}=x_nx_{n-2 \down m + 1} x_{m - 1 \updown 1},
    &n-2\geq m\geq 3,\\
    y^{10}_{m}&{}=x_{n-1}x_nx_{n-2 \down m + 1} x_{m - 1 \updown 1}, &n-2\geq m\geq 3,\\
    y^{11}&{}=x_1 x_{3 \up n - 2} x_{n \down 1} &\\
    y^{12}&{}=x_1 x_{3 \up n - 2} x_{n \down 2},&\\
    y^{13}&{}=x_nx_{n-2 \down 3} x_1,&\\
    y^{14}&{}=x_{n-1}x_nx_{n-2 \down 3} x_1,&\\
    y^{15}&{}=x_{n-2}x_n x_{n - 3 \updown 1},&\\
    y^{16}&{}=x_{n-1}x_{n-2}x_n x_{n - 3 \updown 1},&\\
    y^{17}&{}=x_{n}x_{n-1}x_{n-2}x_n x_{n - 3 \updown 1}.&
  \end{align*}
  Let $x$ be a monomial in $x_1,\dotsc,x_n$ of length $\ell$.
  Then $x$ is a linear combination of monomials $y \in \MD$ with
  $\length (y) \leq \ell$.
\end{theorem}

For convenience in the proof of Theorem~\ref{t| L1 monomials}, we
extend the above definition: we define $y^2_{n-1,m} = y^1_{n,m}$ and
$y^4_{n-1,m} = y^3_{n,m}$, for $m \leq n - 1$; and $ y^3_{2,2}= x_1$,
$y^7_{2}= x_1$, and $y^8_{n-1}= y^7_{n}$.

\begin{proof}
  Let $x$ be a monomial in $x_1,\dotsc,x_n$ of length $\ell$. If $\ell=1$
  then $x\in \MD$. If $\ell=2$, then either $x=0$ or $x\in\MD$ or
  $-x\in\MD$ (using Eqs.~\eqref{e| L} and~\eqref{e| R1}). We use induction
  on $\ell$ and may assume $\ell > 2$.

  We have an $i\in\{1,\dotsc,n\}$ and a monomial $y$ of length $\ell-1$
  such that $x=x_iy$. Moreover, because of the induction hypothesis we
  can assume $y\in\MD$. We will consider $x=x_iy^j_*$ for
  each of the seventeen classes in $\MD$ separately.  In
  each case we will write $x$ as a linear combination of monomials of
  length at most $\ell$, where all monomials in the linear combination of
  length $\ell$ are members of $\MD$. By the induction hypothesis,
  this suffices to prove the theorem.

  We will work modulo monomials of length at most $\ell-1$, so because
  of extremality of $x_j$,
  \begin{equation}
    \label{eq:extr} \tag{XT}
    x_j x_j t = 0
  \end{equation}
  whenever the left hand side occurs in a monomial of length $\ell$.

  We show the cases $j=1$ and $j=3$ as an example, the other cases are
  very similar. For the complete proof, we refer to \cite{cit:postm07}.

  \begin{caselist}
  \item $j=1$, $k\in\{1,\dotsc, n\}$ and $m\in\{1,\dotsc, k\}$. Since
    $\ell>2$, we know that $m<k$.  We distinguish the following
    sub-cases:
    \begin{itemize}
    \item If $i>k+1$ and $k=n-2$, then $x =y^5_m$.
    \item If $i>k+1$ and $k\neq n-2$, then $x = 0$ by Eq.~\eqref{e|
        Q4}.
    \item If $i=k+1$, then $x =y^1_{k+1,m}$.
    \item If $i=k$, then extremality of $x_i$ shows that $x = 0$.
    \item If $i=k-1$, then:
      \[
      x =
      \begin{cases}
        x_{i} x_{i+1} x_{i} \isref{e| L} - x_{i} x_{i}
        x_{i+1} \isref{eq:extr} 0, & \text{if $m=k-1$,} \\
        x_{i} x_{i+1} x_{i}y^1_{i-1,m} \isref{e| P2} 0, &
        \text{otherwise.}
      \end{cases}
      \]
    \item If $i<k-1$ and $i= n-2$, then
      $x =y^2_{n-2,m}.$
    \item If $i<k-1$ and $i\neq n-2$, then $
      x = x_{i} x_{k \down m}.$
      Applying Eq.~\eqref{e| Q1} a sufficient number of times leads to one of the following situations:
      \begin{itemize}
      \item If $i>m$, then
        $x = x_{k \down i+2} x_i x_{i+1} x_{i \down m}\isref{e| P2}0$.
      \item If $i=m$ and $i\neq 1$, then
        $ x =  x_{k \down m+2}x_m x_{m+1} x_m \stackrel{\eqref{e| L},\eqref{eq:extr}}{=} 0$.
      \item If $i=m$ and $i=1$, then
        $x =  x_{k \down 4}x_1 x_{3} x_2 x_1 \isref{e| L}-x_{k \down 4}x_1 x_{3} x_1 x_2\isref{e| P2}0.$
      \item If $i=m-1$ and $i\neq 1$, then
        $x =  x_{k \down m+1} x_{m-1} x_m   \isref{e| L}- y^1_{k,m-1}$.
      \item If $i=m-1$ and $i = 1$, then
        $x = x_{k \down 4}x_1 x_{3} x_2  \isref{e| J} y^3_{k,3}- y^1_{k,1}$.
      \item  If $i<m-1$ and $m \neq 3$, then
        $x = x_{k \down m+1} x_i x_m \isref{e| R1}0$.
      \item  If $i<m-1$ and $m=3$, then
        $x = x_{k \down 4} x_1 x_3 \isref{e| L}-y^7_k.$
      \end{itemize}
    \end{itemize}
 \item $j=3$, $k\in\{3,\dotsc, n\}$ and $m\in\{3,\dotsc, k\}$.
    \begin{itemize}
    \item If $i>k+1$ and $k=n-2$, then $x =y^{9}_m$.
    \item If $i>k+1$ and $k\neq n-2$, then $x = x_{i} x_{k \down m+1} x_{m-1 \updown 1} \isref{e| Q4}
     0.$
    \item If $i=k+1$, then $x =y^{3}_{i,m}$.
    \item If $i=k$, then
      \[
      x = x_{i} y^3_{i,m} =
      \begin{cases}
        x_{i} x_{i-1} x_{i} y^3_{i-1,i-1}\isref{e|
            P2}0, &\text{if $k=m$,} \\
        x_{i} x_{i} y^3_{i-1,m}\isref{eq:extr}
        0,&  \text{otherwise.}
      \end{cases}
      \]
    \item If $i=k-1$, then we have
      \[
      x = x_{i} y^3_{i+1,m} =
      \begin{cases}
        x_{i} x_{i} y^3_{i+1,i}\isref{eq:extr}0, &\text{if $i+1=m$,} \\
        x_{i} x_{i+1} y^3_{i,i}=y^3_{i+1,i+1}, & \text{if $i+1=m+1$,} \\
        x_{i} x_{i+1} x_{i} y^3_{i-1,m}\isref{e| P2}0, & \text{otherwise.}
      \end{cases}
      \]
    \item If $i<k-1$ and $i=n-2$, then $x =y^{4}_{n-2,m}$.
    \item If $i<k-1$ and $i\neq n-2$, then $
      x = x_{i} x_{k \down m+1} x_{m-1 \updown 1}.$
      Repeated application of Eq.~\eqref{e| Q1} leads to one of
      the following situations:
      \begin{itemize}
      \item If $i>m$, then $x = x_{k \down i+2}x_{i} x_{i+1} x_i
        y^3_{i-1,m} \isref{e| P2}0$.
      \item If $i=m$, then
        $ x = x_{k \down m+2} x_{m} x_{m+1} x_{m-1 \updown 1}= y^3_{k,m+1}$.
      \item If $i=m-1$, then $x =x_{k \down m+1} x_{m-1} x_{m-1}
        y^3_{m,m-1}\isref{eq:extr}0$.
      \item If $i<m-1$, $m\neq 3$ and $i\neq 1$, then $x = x_{k \down m+1}
        x_{m-1 \updown i+2}x_{i+3}x_i x_{i+1} x_{i+2}
        x_iy^3_{i+1,i} \isref{e| Q2}0$.
      \item If $i<m -1$, $m\neq 3$ and $i=1$, then
        $ x = x_{k \down m+1} x_{m-1 \updown 4}x_{5}x_1 x_{3}
        x_{4} x_2x_3x_1=0.
        $
        This last identity follows from:
        \begin{equation*}
         \hspace*{1.5cm} x_1 x_{3} x_{4} x_2x_3x_1\isref{e| J} x_1
          x_{3} x_{4} x_1x_3x_2+x_1 x_{3} x_{4} x_3x_2x_1
          {}\stackrel{\eqref{e| Q3},\eqref{e| P2}}{=}{}  x_1 x_{4}
          x_{3} x_1x_3x_2\isref{e| P2}0.\label{e| y3}
        \end{equation*}
      \item If $i<m-1$ and $m=3$, then
        $ x =x_{k \down 4} x_1 x_{2} x_3x_1
        \isref{e| L}-x_{k \down 4} x_1
        x_{2} x_1x_3 \isref{e| P2}0$.
      \end{itemize}
    \end{itemize}
  \end{caselist}
\end{proof}

\subsection{Monomials for $\GB$, $\GA$ and $\GC$}
Recall that in $\FB$, we divide out relations~\eqref{e| R1}
and~\eqref{e| R2} from page~\pageref{e| R1}. In Theorem~\ref{t| L2
  monomials}, we will give a list of $2n^2-3n+1$ monomials in
$x_1,\dotsc,x_n$ and show that these \underline{monomials} span $\FB$ linearly. We
will prove this using Theorem~\ref{t| L1 monomials}.

\begin{theorem}\label{t| L2 monomials}
  Let $\MB$ be the set $\MD$ without the following monomials:
  \begin{enumerate}
  \item $y^1_{n,n-1}= x_n x_{n-1}$,
  \item $y^3_{n,n} = x_{n-1 \updown 1}$,
  \item $y^6_m=x_{n-1}x_nx_{n-2 \down m}$ with $n-2\geq m\geq 1$,
  \item $y^{10}_m=x_{n-1}x_nx_{n-2 \down m+1} x_{m-1 \updown 1}$
    with $n-2\geq m\geq 3$,
  \item $y^{12}=x_{1}x_{3 \up n-2} x_{n \down 3} x_1$,
  \item $y^{14}=x_{n-1}x_nx_{n-2 \down 3} x_1$,
  \item $y^{17}=x_{n}x_{n-1}x_{n-2}x_n x_{n-3 \updown 1}$.
  \end{enumerate} 
  Let $x$ be a monomial in $x_1,\dotsc,x_n$ of length $\ell$.
  Then $x$ is a linear combination of monomials $y \in \MB$ with
  $\length (y) \leq \ell$.
\end{theorem}
\begin{proof}
  Because of Theorem~\ref{t| L1 monomials} it suffices to prove
  that the seven given classes of monomials can be written as a linear
  combination of shorter monomials and monomials from $\MB$.
  We will consider these classes separately and
  work modulo shorter monomials.
  \begin{caselist}
  \item $y = y^1_{n,n - 1} = x_n x_{n-1}\isref{e| R2} 0$.
  \item $y = y^3_{n,n}$. We find $ y = x_{n-1} x_n x_{n-2 \updown 1}
    \stackrel{\eqref{e| Q1}, \eqref{e| R2}}{=} x_n x_{n-1} x_{n-2}
    x_{n-1} x_{n-3 \updown 1} \isref{e| P2} 0$.

  \item $y = y^6_m$. We find $ y = x_{n-1} x_n x_{n-2 \down m}
    \isref{e| J}x_{n} x_{n-1} y^1_{n-2,m} + (x_{n-1} x_n) y^1_{n-2,m}
    \isref{e| R2}y^1_{n,m}$.

  \item $y = y^{10}_m$. We find $ y = x_{n-1} x_n x_{n-2 \down m+1} x_{m-1 \updown 1} \isref{e| J}x_{n} x_{n-1} y^3_{n-2,m}
    + (x_{n-1} x_n) y^3_{n-2,m} \isref{e| R2}y^3_{n,m}$.

  \item $y = y^{12}$. $y^{12}$ can be
    written as a linear combination of the following monomials:
    \[
    y^2_{3,1}, y^3_{n,n-1}, y^4_{3,3}, y^4_{4, 4}, \dotsc,
    y^4_{n-2,n-2}, y^8_{2}, y^{12}, y^{16} \text{ and } x_{n-2} x_{n-1}
    x_n x_{n-3 \updown 1} \eqqcolon t.
    \]
    All of these except $t$ are in $\MB$, so we only need to
    analyze $t$.
    \[
    t = x_{n-2} x_{n-1} x_n x_{n-3 \updown 1} \isref{e|
      Q1}x_{n-2}x_{n}x_{n-1} x_{n-3 \updown 1} = y^4_{n-2,n-2}.
    \]

  \item $y = y^{14}$. We find $ y = x_{n-1} x_n x_{n-2 \down 3} x_1
    \isref{e| Q1} x_{n \down 3} x_1 = y^7_n$.

  \item $y = y^{17}$. We find $ y = x_n x_{n-1} x_{n-2} x_n x_{n-3
      \updown 1} \isref{e| Q1} x_{n-1} x_n x_{n-2} x_n x_{n-3 \updown
      1} \isref{e| P2} 0$. \qedhere
  \end{caselist}
\end{proof}

The following two theorems can be proved in exactly the same manner.

\begin{theorem} \label{t| L3 monomials}
  Let $\MA$ be the set consisting of the following
  monomials:
  \begin{align*}
    y^1_{k,m}&{}=x_{k \down m},&n\geq k\geq m\geq 1,\\
    y^3_{k,m}&{}=x_{k \down m+1}x_{m-1 \updown 1}, & n \geq k\geq
    m\geq 3, \\
    y^7_{k}&{}= x_{k \down 3} x_1,&n\geq k\geq 3.
  \end{align*}
  Let $x$ be a monomial in $x_1,\dotsc,x_n$ of length $\ell$. Then
  $x$ is a linear combination of monomials $y \in \MA$ with $\length
  (y) \leq \ell$.
\end{theorem}

\begin{theorem}\label{t| L4 monomials}
  Let
  \[
  \MC = \{y^1_{k,m} = x_{k \down m} \mid n \geq k \geq m \geq 1 \}.
  \]
  Let $x$ be a monomial in $x_1, \dotsc, x_n$ of length $\ell$. Then
  $x$ is a linear combination of monomials $y \in \MC$ with $\length
  (y) \leq \ell$.
\end{theorem}

\section{The parameter space}\label{s:parameter-space}
Recall that $\Xx{\Gamma} = \{\ff \mid \dim \Ll{\Gamma}{\ff} =
\card{\MM_\Gamma}\}$. For $\Gamma \in \{\GA, \GB, \GC, \GD\}$, we will
find bijections $\psiG$ from $\Xx{\Gamma}$ to a vector space, such that
$\Ll{\Gamma}{\ff}$ is isomorphic to the Lie algebra of the
corresponding Chevalley
type if $\psiG (\ff)$ is in a certain open dense subset of that vector
space.  In this section, we will be constructing the bijection. It
will take some work to show that the map is injective; this will be
the content of Lemmas~\ref{lm:f-D}, \ref{lm:f-B}, \ref{lm:f-A}
and~\ref{lm:f-C}.  We will find the open dense subset in
Section~\ref{s:alwaysours}.

\subsection{Parameters for $\GD$}
\begin{lemma}\label{lm:f-D}
  Let
  \begin{multline*}
    \psiD \colon \XD \to \FF^{n+4}, \ff \mapsto \\ (\ff_1 (x_2), \ff_2
    (x_3), \dotsc, \ff_{n-1} (x_n); \ \ff_1 (x_3), \ff_1 (x_2 x_3),
    \ff_{n-2} (x_n), \ff_{n-2} (x_{n-1} x_n), \ff_1 (x_{3 \up n-2}
    x_{n \down 2})).
  \end{multline*}
  Then $\psiD$ is injective.
\end{lemma}
\begin{proof}
  The values of all $\ff_i$ together determine the values of the
  extremal bilinear form on all of $\FD$, since $f (x_i y, z) = f
  (x_i, y z) = \ff_i (y z)$. We will show that each value $\ff_i (y)$
  for $y \in \MD$ can be expressed in the values $\ff_j (z)$ in the
  theorem. To make this notationally convenient, let $\FFD$ be the
  rational function field obtained from $\FF$ by extending it with $n
  + 4$ symbols as follows. For every $\ff_j (z)$ in the lemma, we
  extend $\FF$ with the symbol $\ff_j (z)$ and assume that evaluating
  $\ff_j$ at $z$ yields this symbol $\ff_j (z) \in \FFD$. We will show
  that each value $\ff_i (y)$ is an element of $\FFD$ for \emph{all} $i$ and
  \emph{all} $y \in \MD$.

  Let $y = y^j_* \in \MD$ of length $\ell$ and $i \in \{1, \dotsc,
  n\}$.  We will use induction on $\ell$.  We consider the seventeen
  classes of monomials in $\MD$ separately, and give cases $j=1$ and $j=3$
  as an example. For the complete proof, we refer to \cite{cit:postm07}.

   Let $x = x_i x_i y = f(x_i, y) x_i$. It is sufficient to prove that $x=0$ or $f(x_i,y)=0$.
  \begin{caselist}
  \item $j=1$, $k\in\{1,\dotsc,n\}$, $m\in\{1,\dotsc,k\}$. Note that
    $y_{k,m}^1 = x_{k \down m}$.

    If $m=k$ then $y$ has length $1$ and $f (x_i, y)$ is either $0$ or in
    the list of values in the theorem. So assume that $k>m$.
    \begin{itemize}
    \item If $i>k+1$ and $k=n-2$, then:
      \begin{align*}
        x \isref{e| L}&- x_{n}x_nx_{n-2 \down m+2}
        x_{m} x_{m+1}
        \isref{e| Q1} - x_{m}x_{n}x_nx_{n-2 \down m+1} \\
        \isref{eq:extr} & -f (x_n, y^1_{n-2,m+1})x_mx_n\isref{e| R1}0.
      \end{align*}

    \item If $i>k+1$ and $k < n-2$, then
      \[
      x \isref{e| Q1}x_{k}x_ix_{i}y^1_{k-1,m}
      \isref{eq:extr}  f (x_i, y^1_{k-1,m})x_kx_i\isref{e| R1} 0.
      \]

    \item If $i=k+1$, $m=k-1$ and $i=n$, then
      $f (x_i, y)=f (x_n, x_{n-1}x_{n-2})\in\FFD$.

    \item If $i=k+1$, $m=k-1$ and $3<i<n$, then
      \[
      f (x_i, y)=f(x_i,x_{i-1}x_{i-2})
      \isref{e| f} - f(x_{i-1},x_{i}x_{i-2})\isref{e| R1}0.
      \]

    \item If $i=k+1$, $m=k-1$ and $i=3$ then $f (x_i, y)=f (x_3,x_{2}x_{1})\in\FFD$.

    \item If $i=k+1$, $m<k-1$ and $m>1$, then
      \[
      x  \isref{e| L} - x_i x_{i \down m + 2} x_m x_{m+1}
      \isref{e| Q1} - x_m x_i x_{i \down m + 1} \isref{eq:extr} -
      f (x_i, y^1_{i-1,m+1}) x_m x_i \isref{e| R1} 0.
      \]

    \item If $i=k+1$, $m < k-1$ and $m=1$, then
      \begin{align*}
        x \isref{e| J}& x_i x_{i \down 4} x_2 x_3 x_1 + x_i x_{i \down 4} x_1 x_2 x_3 \isref{e| Q1} x_2 x_i x_{i \down 3} x_1 + x_1
        x_i x_{i \down 4} x_2 x_3 \\
        \isref{eq:extr} & f (x_i, x_{i-1 \down 3} x_1) x_2 x_i + f (x_i, x_{i-1 \down 4} x_2 x_3) x_1 x_i
        \isref{e| R1} 0.
      \end{align*}

    \item If $i=k$, then $x = x_i x_i x_{i \down m} = 0$.

    \item If $i=k-1$ and $i>m$, then
      \[
      f (x_i, y)=f (x_i, x_{i+1}x_{i}y^1_{i-1,m})\isref{e| P5}-f (x_i, x_{i+1})
      f (x_i, y^1_{i-1, m}) \in \FFD.
      \]
      Here we use the induction hypothesis for $f (x_i, y^1_{i-1, m})$.

    \item If $i=k-1$ and $i=m$, then
      \[x = x_i x_i x_{i+1} x_i \isref{e| L} - x_i x_i x_i x_{i + 1} =
      0. \]

    \item If $i<k-1$, $i=n-2$ and $m=n-1$, then
      \[
      f (x_i, y) = f (x_{n-2}, x_n x_{n-1}) \stackrel{\eqref{e| f},
        \eqref{e| ff}}{=} f(x_{n},x_{n-1}x_{n-2})\in\FFD.
      \]

    \item If $i<k-1$, $i=n-2$ and $m=n-2$, then
      \[\hspace*{.75cm}
      f (x_i, y) = f (x_{n-2}, x_n x_{n-1} x_{n-2}) \isref{e| L} - f
      (x_{n-2}, x_n x_{n-2} x_{n-1}) \isref{e| P5} f (x_{n-2}, x_n)
      f (x_{n-2}, x_{n - 1}).
      \]

    \item If $i<k-1$, $i=n-2$ and $m\leq n-3$, then
      \begin{align*}
        f (x_i, y) ={}& f (x_{n-2}, x_{n \down m}) \isref{e|
          f}f(x_{n-2}x_n,x_{n-1 \down m})\\
        \isref{e| L}&-
        f(x_nx_{n-2},x_{n-1 \down m}) \isref{e|
          f}-f(x_n,x_{n-2}x_{n-1}x_{n-2 \down m})\inref{e| P2}\FFD.
      \end{align*}
    \item If $i<k-1$, $i<n-2$ and  $k>3$, then
      \[
      x = x_i x_i x_{k \down m} \isref{e| Q1} x_k x_i x_i x_{k-1 \down m} \isref{eq:extr} f (x_i, y_{k-1,m}^1) x_k x_i \isref{e|
        R1} 0.
      \]

    \item If $i < k-1$, $i<n-2$, $k=3$ and $m = 2$, then $ f (x_i, y) = f (x_1, x_3
      x_2) \in \FFD$.

    \item If $i<k-1$, $i< n-2$, $k=3$ and $m=1$, then we have:
      \[\hspace*{.75cm}
      f (x_i, y) = f(x_{1},x_3x_2x_1)\isref{e|
        f}f(x_{1}x_3,x_2x_1)\isref{e| L} f(x_{3}x_1,x_1x_2) \isref{e|
        f}
      f(x_{3},x_1x_1x_2) \in \FFD.
      \]
    \end{itemize}

  \item $j=3$, $k\in\{3,\dotsc,n\}$, $m\in\{3,\dotsc,k\}$. Remember
    that $y^3_{k, m} = x_{k \down m + 1} x_{m - 1 \updown 1}$.
    \begin{itemize}
    \item If $i>k+1$ and $m=k$, then
      \[
      x=x_{i}x_ix_{k-1}y^3_{k,k-1} \isref{e|
        Q1}x_{k-1}x_ix_{i}y^3_{k,k-1} \isref{eq:extr}
      f (x_{i}, y^3_{k,k-1})x_{k-1}x_i \isref{e| R1}0.
      \]

    \item If $i>k+1$ and $m<k$, then
      \begin{align*}\hspace*{.75cm}
        x={}&x_{i}x_ix_{k \down m+1}x_{m-1}x_m x_{m-2 \updown 1}
        \isref{e| Q1} x_{m-1}x_{i}x_ix_{k \down m+1}x_m x_{m-2 \updown 1}
        \in \FFD x_{m-1}x_i\isref{e| R1}\{0\}.
      \end{align*}

    \item If $ i=k+1$, $m=k$ and $i=n$, then
      \begin{align*}\hspace*{.5cm}
        f (x_i, y) &\ \ \ =\ \ \ f (x_n, x_{n-2 \updown 1}) \isref{e| f}
        f (x_n x_{n-2}, x_{n-1} x_{n-3 \updown 1})
        \isref{e| L}
        - f (x_{n-2} x_n, x_{n-1} x_{n-3 \updown 1}) \\
         &\ \ \isref{e| f}\ \
        - f (x_{n-2}, x_n x_{n-1} x_{n-3 \updown 1})
        \isref{e| Q1}
        - f (x_{n-2}, x_{n-3} x_n x_{n-1} x_{n-2} x_{n-4 \updown 1})\\
        &\hspace*{-.1cm}\stackrel{\eqref{e| f}, \eqref{e| L}}{=}
        f (x_{n-3}, x_{n-2} x_n x_{n-1} x_{n-2} x_{n-4 \updown 1}). \\
        \intertext{Now  we repeat the steps in the last line $n - 6$
          times to obtain}
        f (x_i, y)  &\ \ \ =\ \ \  (-1)^n f (x_3, x_{4 \up n-2} x_{n \down 4} x_2
        x_3 x_1) \\
        &\ \ \ \isref{e| J}\ \ \  (-1)^n (f (x_3, x_{4 \up n-2} x_{n \down 4} (x_2 x_3) x_1) + f (x_3, x_{4 \up n-2} x_{n \down 4} x_3 x_2 x_1)).
      \end{align*}
      We will treat these terms separately.
      \begin{align*}\hspace*{.5cm}
        f (x_3, x_{4 \up n-2} x_{n \down 4} (x_2 x_3) x_1)&
        \ \ \ \isref{e| L}\ \ \ f (x_3, x_{4 \up n-2} x_{n \down 4} x_1
        x_3 x_2)
        \isref{e| Q1} f (x_3, x_1 x_{4 \up n-2} x_{n \down 2}) \\ &\stackrel{\eqref{e| f}, \eqref{e| L}}{=}{}
        -f (x_1, x_{3 \up n-2} x_{n \down 2}) \in \FFD,
        \end{align*}
        and
         \begin{align*}
         f (x_3, x_{4 \up n-2} x_{n \down 1})&
        \stackrel{\eqref{e| f}, \eqref{e| L}}{=}{}
        - f (x_4, x_3 x_{5 \up n-2} x_{n \down 1})
        \isref{e| Q1} - f (x_4, x_{5 \up n - 2} x_{n \down 5} x_3
        x_{4 \down 1}) \\&\ \ \ \inref{e| P1}\ \ \ \FFD.
      \end{align*}

    \item If $i=k+1$, $m=k$ and $i<n$ then
      \[
        x=x_{i}x_ix_{i-2}y^3_{i-1,i-2}
        \isref{e| Q1}x_{i-2}x_ix_{i}y^3_{i-1,i-2}
        \inref{eq:extr} \FFD x_{i-2}x_i\isref{e| R1}0.
      \]

    \item If $i=k+1$ and $m<k$, then
      \[
      x=x_{i}x_{i \down m+1}x_{m-1}y^3_{m,m-1}
      \isref{e| Q1}
      x_{m-1}x_{i}x_{i \down m+1}y^3_{m,m-1} \inref{eq:extr}
      \FFD x_{m-1}x_i\isref{e| R1}0.
      \]

    \item If $i=k$ and $m=k$, then
      \[
      f (x_i, y)=f(x_i,x_{i-1}x_iy^3_{i-1,i-1})\inref{e| P5}\FFD.
      \]

    \item If $i=k$ and $m<k$, then
      \[
      f (x_i, y)=f(x_i,x_iy^3_{i-1,m})\isref{e| f}f(x_ix_i,y^3_{i-1,m})\isref{e| L}0.
      \]

    \item If $i=k-1$ and $m=k$, then
      \[
      f (x_i, y)=f(x_i,x_iy^3_{i+1,i})\isref{e| f}f(x_ix_i,y^3_{i+1,i})\isref{e| L}0.
      \]

    \item If $i=k-1 $ and $m=k-1$, then
      \begin{align*}
        f (x_i, y)={}&f(x_i,x_{i+1}x_{i-1}x_ix_{i-2 \updown 1})\isref{e|
          f}f(x_ix_{i+1},x_{i-1}x_ix_{i-2 \updown 1})\\
        \isref{e| L}&-f(x_{i+1}x_{i},x_{i-1}x_ix_{i-2 \updown 1})\isref{e|
          f}-f(x_{i+1},x_ix_{i-1}x_ix_{i-2 \updown 1}) \inref{e| P2}
        \FFD
      \end{align*}

    \item If $i=k-1$ and $m<k-1$, then
      \[
      f (x_i, y)=f(x_i,x_{i+1}x_iy^3_{i-1,m})\inref{e| P5}\FFD.
      \]

    \item If $i=k-2$, $m=k$ and $i>1$, then
      \[\hspace*{.5cm}
      f (x_i, y) = f (x_i, x_{i + 1} x_{i + 2} x_i x_{i + 1} x_{i-1 \updown 1})
      \inref{e| Q3a} \FFD + f (x_i, x_{i+1} x_i x_{i+2} x_{i+1} x_{i-1 \updown 1}) \isref{e| P5} \FFD.
      \]
    \item If $i=k-2$, $m=k$ and $i=1$, then
      \[
      f (x_1,y) = f (x_1, x_2 x_3 x_1) \isref{e|
        L}-f(x_1,x_{2}x_{1}x_3)\inref{e| P5}\FFD.
      \]

    \item If $i=k-3, $ 
    $m=k$ and $i>1$, then
      \[
      x=x_{i}x_ix_{k-1}y^3_{k,k-1} \isref{e|
        Q1}x_{k-1}x_{i}x_{i}y^3_{k,k-1} \isref{eq:extr}
      f (x_{i}, y^3_{k,k-1})x_{k-1}x_{i}\isref{e| R1}0.
      \]

    \item If $i=k-3$, $m=k$ and $i=1$, then
      \begin{align*}
        f (x_i, y)={}&f(x_1,x_{3}x_4x_2x_3x_1)
        \isref{e| J}f(x_1,x_{3}x_4x_3x_2x_1)+f(x_1,x_{3}x_4x_1x_3x_2)\\
        \isref{e| Q1}&f(x_1,x_{3}x_4x_3x_2x_1)+ f(x_1,
        x_{3}x_1x_4x_3x_2) \stackrel{\eqref{e| P2}, \eqref{e|
            P5}}{\in} \FFD.
      \end{align*}

    \item If $i < k - 3$ and $m=k$, then
      \[
      x=x_{i}x_ix_{k-1}y^3_{k,k-1} \isref{e|
        Q1}x_{k-1}x_{i}x_{i}y^3_{k,k-1} \isref{eq:extr}
      f (x_{i}, y^3_{k,k-1})x_{k-1}x_{i}\isref{e| R1}0.
      \]

    \item If $i<k-1$, $m<k$, $i=n-2$ and $m=n-1$, then
      \[
      f (x_i, y) = f (x_{n-2}, x_n x_{n-2} y^3_{n-1,n-2}) \inref{e| P5}
      \FFD.
      \]

    \item If $i<k-1$, $m<k$, $i=n-2$ and $m=n-2$, then
      \begin{align*}\hspace*{.5cm}
      f (x_i, y) &= f (x_{n-2}, x_n x_{n-1} x_{n-3 \updown 1})
      \stackrel{\eqref{e| f}, \eqref{e| L}}{=}
      -f (x_n, x_{n-2} x_{n-1} x_{n-3 \updown 1}) = - f(x_n,
      y^3_{n-1,n-1})\\
      & \in \FFD,
      \end{align*}
      as proven earlier.

    \item If $i<k-1$, $m<k$, $i=n-2$ and $m<n-2$, then
      \begin{align*}
        x&={}x_{n-2}x_{n-2}x_{n \down m+1}x_{m-1}y^3_{m,m-1} \isref{e|
          Q1} x_{m-1}x_{n-2}x_{n-2}x_{n \down m+1}y^3_{m,m-1} \\
       &\hspace*{-.2cm}\isref{eq:extr} f (x_{n-2}, y^3_{n,m-1})x_{m-1}x_{n-2}
        \isref{e| R1} 0.
  \end{align*}
    \item If $i<k-1$, $m<k$ and $i< n-2$, then
      \[
      x = x_i x_i x_k y^3_{k-1,m} \isref{e| Q1} x_k x_i x_i
      y^3_{k-1,m} \isref{eq:extr} f (x_i, y^3_{k-1,m}) x_k x_i \isref{e|
        R1} 0.
      \]

    \end{itemize}
  \end{caselist}
\end{proof}
\subsection{Parameters for $\GB$}
\begin{lemma}\label{lm:f-B}
  Let
  \begin{multline*}
    \psiB\colon \XB \to \FF^{n+2}, \ff \mapsto \\ (\ff_1 (x_2), \ff_2
    (x_3), \dotsc, \ff_{n-2} (x_{n - 1}); \ \ff_1 (x_3), \ff_1 (x_2
    x_3), \ff_{n-2} (x_n), \ff_1 (x_{3 \up n-2} x_{n \down 2})).
  \end{multline*}
  Then $\psiB$ is injective.
\end{lemma}
\begin{proof}
  Remember that $\FB$ is a quotient of $\FD$. Hence relations between
  values of $\ff_i$ that hold in $\FD$, hold in $\FB$ as well. This
  allows us to express $\ff_i (y)$ in the values of
  Lemma~\ref{lm:f-D}. It then suffices to prove that $f (x_n,
  x_{n-1})$ and $f (x_n, x_{n-1} x_{n-2})$ are zero:
  \begin{itemize}
  \item $f (x_n, x_{n-1})=0$, because $x_nx_nx_{n-1}\isref{e| R2}0$,
  \item $f (x_n, x_{n-1}x_{n-2})={}f(x_n,x_{n-1}x_{n-2})\isref{e|
      f}f(x_nx_{n-1},x_{n-2})\isref{e| R2}0.$ \qedhere
  \end{itemize}

\end{proof}
\subsection{Parameters for $\GA$}
\begin{lemma}\label{lm:f-A}
  Let
  \[
  \psiA\colon \XA \to \FF^{n+1}, \ff \mapsto (\ff_1 (x_2), \ff_2 (x_3),
  \dotsc, \ff_{n-1} (x_n); \ \ff_1 (x_3), \ff_1 (x_2 x_3)).
  \]
  Then $\psiA$ is injective.
\end{lemma}
\begin{proof}
$\FA$ is a quotient of $\FD$, so by Lemma~\ref{lm:f-D} and a
  reasoning similar to the one in the proof of Lemma~\ref{lm:f-B},
it suffices to show that $f (x_n, x_{n-2})$, $f
  (x_n, x_{n-1} x_{n-2})$ and $f (x_1, x_{3 \up n-2} x_{n \down 2})$
  are zero:
  \begin{itemize}
  \item$f (x_n, x_{n-2})=0$, because $x_n x_n x_{n-2} \isref{e| R3} 0$,
  \item$f (x_n, x_{n-1}x_{n-2}) \isref{e| L} - f (x_n, x_{n-2}
    x_{n-1}) \isref{e| f} - f (x_n x_{n-2}, x_{n-1}) \isref{e| R3} 0$,
  \item $f (x_1, x_{3 \up n-2} x_{n \down 2}) \isref{e| f} f (x_1x_3, y^2_{4,2})\isref{e| R4} 0$. \qedhere
  \end{itemize}
\end{proof}
\subsection{Parameters for $\GC$}
\begin{lemma}\label{lm:f-C}
  Let
  \[
  \psiC\colon \XC \to \FF^{n-1}, \ff \mapsto (\ff_1 (x_2), \ff_2 (x_3),
  \dotsc, \ff_{n-1} (x_n)).
  \]
  Then $\psiC$ is injective.
\end{lemma}
\begin{proof}
  $\FC$ is a quotient of $\FA$, so by Lemma~\ref{lm:f-A} and a
  reasoning similar to the one in the proof of Lemma~\ref{lm:f-B}, it
  suffices to show that $f (x_3, x_1)$ and $f (x_3, x_2 x_1)$ are
  zero:
  \begin{itemize}
  \item$f (x_3,x_{1})=0$, because $x_3 x_3 x_1 \isref{e| R4} 0$,
  \item$f (x_3,x_{2}x_{1})\isref{e| L} - f (x_3, x_1 x_2) \isref{e| f}
    - f (x_3 x_1, x_2) \isref{e| R4} 0$. \qedhere
  \end{itemize}
\end{proof}

\begin{corollary}
  $\psiG (\Xx{\Gamma})$ is an algebraic variety for all $\Gamma \in
  \{\GA, \GB, \GC, \GD\}$.
\end{corollary}

\section{Realizations of the four classical families}\label{s:realizations}
In this section, we will find generating sets of extremal elements for
the four classical families of Lie algebras, where these generating
sets realize the graphs $\GA$, $\GB$, $\GC$ and $\GD$.  We keep $n$ as
the number of extremal generators and will see that these graphs
correspond to the Lie algebras of type~$\Aa$, $\Bb$, $\Cc$ and~$\Dd$,
respectively.  In particular, the objective of this section will be
the formulating and proving of Theorems~\ref{thm:realizeslinear},
\ref{thm:realizestwt}, \ref{thm:realizestwotr}
and~\ref{thm:realizestrspt}, giving these explicit generators.

The extremal elements in Lie algebras of type $\Aa$ or $\Cc$
correspond to \emph{transvections}, which will be discussed in
Section~\ref{ss:transvections}. In the orthogonal Lie algebras, the
extremal elements correspond to the \emph{Siegel
  transvections}\index{Siegel transvection} or Siegel
transformations\index{Siegel transformation|see{Siegel transvection}}.
These are examined in Section~\ref{ss:siegel-transvections}. In both
sections, we first explore these elements in a general setting, and
then discuss generators for the two series of Lie algebras
specifically.

\subsection{Transvections}\label{ss:transvections}
Let $n \in \None$, $x \in V = \FF^n$, $h \in V^*$, and fix a basis
$e_1, \ldots, e_n$ of $V$ and a corresponding dual basis 
$f_1, \ldots, f_n$ of $V^*$.
We will see that the linear
transformation $x \tensor h\colon v \mapsto h (v)\ x$ is an extremal
element of $\sll (V)$ if $h (x) = 0$ and $x, h$ nonzero. A
\emph{transvection}\index{transvection} is a linear transformation of the form $1 + x
\tensor h$ where $h (x) = 0$ and $x, h$ nonzero. We call $x$ the
\emph{centre}\index{centre} of the transvection and $h$ the \emph{axis}\index{axis}. We then
call $x \tensor h$ an \emph{infinitesimal transvection}\index{infinitesimal transvection}.

A \emph{transvection group}\index{transvection group} is a group $\{1
+ t x \tensor h \mid t \in \FF\}$. The Lie algebra of a transvection
group consists of the transvections $t x \tensor h$. We will use a
result of McLaughlin~\cite{cit:mclau67} which classifies groups
generated by transvection subgroups. This is a weaker version of a
reformulation by Cameron and Hall, Theorem~2 from~\cite{cit:camha91}:
\begin{theorem} \label{thm:camhall}
  Let $G$ be a nontrivial group of linear transformations of the
  finite-dimensional $\FF$-vector space $V$, which is generated
  by $\FF$-transvection subgroups. If $V$ is spanned by a
  $G$-orbit on centres of these transvection subgroups, and $V^*$ is
  spanned by the axes, then one of the following holds:
  \begin{enumerate}
  \item $G = \SL (V)$;
  \item $G = \Sp (V, B)$ for some symplectic form $B$.
  \end{enumerate}
\end{theorem}
We will need a tool to distinguish between $\SL (V)$ and $\Sp (V,
B)$. This tool will be provided by analysis of the occurrence of
Heisenberg subalgebras generated by pairs of extremal elements,
further detailed in Lemma~\ref{lm:classify-transvection-generables}
and Corollary~\ref{cor:sp-no-exspp}.

The following lemmas are easily seen to be true.
\begin{lemma} \label{lm:ext-spans-sl}
  $\sll(V)$ is the span of the infinitesimal transvections in $\gll(V)$.
\end{lemma}
\begin{lemma} \label{lm:inftrans-are-ext}
  Infinitesimal transvections are extremal elements of $\sll (V)$.
\end{lemma}
\begin{lemma} \label{lm:all-extr-are-inf-tr-sl}
  All extremal elements in $\sll (V)$ are infinitesimal transvections.
\end{lemma}
\begin{lemma} \label{lm:classify-transvection-generables}
  Let $ x \tensor h$ and $ y \tensor k$ be two infinitesimal
  transvections and let $\LL = \langle x \tensor h, y \tensor k
  \ranLie$. Then the isomorphism class of $\LL$ depends on the
  geometrical configuration of $\FF x$, $\FF y$, $\Ker h$ and $\Ker
  k$, as follows:
  \begin{itemize}
  \item If $\FF x = \FF y$ and $\Ker h = \Ker k$, then $\LL$ is
    one-dimensional.
  \item If either $\FF x = \FF y$ or $\Ker h = \Ker k$ but not both,
    then $\LL$ is two-dimensional.

    \noindent Assume for the other cases that $\FF x \not= \FF y$ and $\Ker h
    \not= \Ker k$.
  \item If $\FF x \subset \Ker k$ and $\FF y \subset \Ker h$, then
    $\LL$ is two-dimensional as in the preceding case.
  \item If either $\FF x \subset \Ker k$ or $\FF y \subset \Ker h$ but
    not both, then $\LL$ is isomorphic to the Heisenberg algebra.
  \item If  $\FF x \not\subset \Ker k$ and $\FF y \not\subset \Ker h$,
    then $\LL$ is isomorphic to $\sll_2$.
  \end{itemize}
\end{lemma}
If $x \tensor h$ and $y \tensor k$ generate a Heisenberg algebra, we
say that they form an \emph{extraspecial pair}\index{extraspecial pair}.

We will also realize $\spp (V)$, the Lie algebra of type $\Cc$, using
infinitesimal transvections. In order to do this, let us assume that
$n$ is even and that we have a nondegenerate symplectic form $B$. We
will denote the matrix of $B$ by $B$ as well. For $y \in V$, we write
\[
u (y) \coloneqq y \tensor (v \mapsto B (y, v)) \colon V \to V, v
\mapsto B (y, v) y.
\]
Hence $u (y)$ is an infinitesimal transvection. The following lemma is
easily proven:
\begin{lemma}
  The infinitesimal transvections in the Lie algebra of type $\Cc$ can
  all be written as $u (y)$ for some $y$.
\end{lemma}

Since $u (y)$ is an extremal element in $\sll (V)$, it is also an
extremal element in $\spp (V)$. The following two lemmas are
equivalents of Lemmas~\ref{lm:ext-spans-sl}
and~\ref{lm:all-extr-are-inf-tr-sl} for $\spp (V)$.

\begin{lemma} \label{lm:ext-spans-sp}
  $\spp (V)$ is spanned by its infinitesimal transvections.
\end{lemma}
\begin{lemma}
  All extremal elements in $\spp (V)$ are infinitesimal transvections.
\end{lemma}
\begin{proof}[Proof sketch]
  We may assume that $B = \left(\begin{smallmatrix} 0 & I \\ -I & 0
    \end{smallmatrix} \right)$. Let $M$ be the matrix of an extremal
  element of $\spp (V)$, then $M$ can be written as a block matrix
  $\left(
    \begin{smallmatrix}
      M_{11} & M_{12} \\ M_{21} & -M_{11}^T
    \end{smallmatrix}
  \right)$, where $M_{12}$ and $M_{21}$ are symmetric. For every
  matrix $A$ also in $\spp (V)$, we have that
  \begin{equation}
    [M, [M, A]] = M^2 A - 2 MAM + AM^2 \in \FF M.\label{eq:extremal-matrix}
  \end{equation}
  We will mostly take for $A$ a block matrix $\left(
    \begin{smallmatrix}
      A_{11} & A_{12} \\ A_{21} & A_{22}
    \end{smallmatrix}
  \right)$ with $A_{12}$ and $A_{21}$ symmetric and $A_{22} = -
  A_{11}^T$. However, substituting the $2n\times 2n$ identity matrix
  into the linear equation~\eqref{eq:extremal-matrix} yields a
  tautology, so we can add multiples of the $n\times n$ identity
  matrices to $A_{11}$ and $A_{22}$ independently. We take the
  following values for $A$.
  \begin{itemize}
  \item Take $A_{11} = A_{12} = A_{22} = 0$ and for $A_{21}$ all
    matrices $E_{ij}$ that has a one in position $(i, j)$ and zeroes
    elsewhere. This shows that we can find a vector $y_0 \in \FF^n$
    such that $M_{12} = y_0 y_0^T$. Similarly there is an $y_1 \in \FF^n$
    with $M_{21} = - y_1 y_1^T$.
  \item Take $A_{11} = I$ and $A_{12} = A_{21} = A_{22} = 0$, then we
    find that
    \begin{equation}\hspace*{.5cm}
      M' \coloneqq \begin{pmatrix}
        2 M_{12} M_{21} & - M_{11} M_{12} - M_{12} M_{11}^T \\
        - M_{21} M_{11} - M_{11}^T M_{21} & - 2 M_{21} M_{12}
      \end{pmatrix} \qquad \text{is a multiple of
        $M$.}\label{eq:matrix-for-m11}
    \end{equation}
    If $M'$ is nonzero, it follows fairly easily that $M_{11} = - y_0
    y_1^T$, and we are done. We continue with the case $M' = 0$, where
    one sees fairly easily that the products of pairs of distinct
    elements of $\{M_{11}, M_{12}, M_{21}\}$ are zero.
  \item Take $A_{12} = A_{21} = 0$ and $A_{22} = - A_{11}^T$. Using
    the fact that the products mentioned above are zero, we obtain
    that
    \begin{multline}\hspace*{1cm}M''\coloneqq
      \begin{pmatrix}
        [M_{11}, [M_{11}, A_{11}]] + 2 M_{12} A_{11}^T M_{21} &
        -2 (M_{11} A_{11} M_{12} + M_{12} A_{11}^T M_{11}^T) \\
        -2(M_{21} A_{11} M_{11} + M_{11}^T A_{11}^T M_{21}) &
        - [M_{11}^T, [M_{11}^T, A_{11}^T]] - 2 M_{21} A_{11} M_{12}
      \end{pmatrix}  \\
      \quad\text{is a multiple of $M$.}\label{eq:mat2-for-m11}
    \end{multline}
    We take for $A_{11}$ the matrices $E_{ij}$. The rest of the proof
    is easy if one distinguishes the case where $y_0 = y_1 = 0$ from
    the other case.\qedhere
  \end{itemize}
\end{proof}
\begin{corollary} \label{cor:sp-no-exspp}
  $\spp (V)$ does not contain an extraspecial pair.
\end{corollary}
\begin{proof}
  Let $x \tensor h$ and $y \tensor k$ form an extraspecial pair. Then
  we may assume that $h (y) = 0 \not= k (x)$. But $h (y) = B (x, y) =
  - B (y, x) = - k (x)$.
\end{proof}

We now proceed with the theorems giving the generating extremal
elements. 
\begin{theorem} \label{thm:realizeslinear}
  Suppose that $n$ is even. Let $B$ be the nondegenerate symplectic
  form determined by the matrix $\left(\begin{smallmatrix} 0 & I \\ -I
      & 0 \end{smallmatrix}\right)$. The transformations
  \[ u (x_i) \colon v \mapsto B (x_i, v) x_i, \qquad i = 1, \dotsc, n, \]
  realize the graph $\GC$ for $n \geq 2$
  if we take these values for $x_i$:
  \begin{align*}
    x_{2 \ell - 1} &= e_\ell & \text{for $1 \leq \ell \leq n/2$,} \\
    x_{2 \ell} &= e_{\ell + n/2} + e_{\ell + n/2 + 1} & \text{for $1 \leq
      \ell < n/2$,} \\
    x_n &= e_n.
  \end{align*}
  The Lie algebra $\LL = \langle u (x_i) \ranLie$ is $\spp_n$.
\end{theorem}
\begin{proof}
  It is easy to check that the transformations $u (x_i)$ realize the
  graph.  In order to prove that they generate $\spp_n$, consider the
  group $G = \langle 1 + t u (x_i) \mid t \in \FF, 1 \leq i \leq n
  \ranGp$, of which $\LL$ is the Lie algebra. The action of $G$ on
  $\LL$ is such that
  \[
  u (x)^{1 + t u (y)} = (1 + t u (y)) u (x) (1 - t u (y)) = u (x^{1 +
    t u (y)}),
  \]
  as can be seen immediately by inspection. So if $B (x, y) \not= 0$,
  then
  \[
  (1 + B (x, y)^{-1} u (y)) (1 + B (x, y)^{-1} u (x)) (y) = x.
  \]
  This shows that the orbit of $G$ on $x_1$ spans $V$. Hence, using
  Theorem~\ref{thm:camhall}, either $G = \Sp (V, B)$ or $G = \SL (V)$.
  Thus $\LL$ is either $\spp_n$ or $\sll_n$. But all given
  transformations are in $\spp_n$. This can be verified by examining
  the matrices $A_i^T M + M A_i$, where $M$ is the matrix of $B$ and
  $A_i$ is the matrix of $u (x_i)$: if we view elements of $V$ as
  column vectors, then $A_i = x_i x_i^T M$, so $A_i^T = - M x_i
  x_i^T$.  Thus $A_i^T M + MA_i = 0$, so $\LL = \spp_n$.
\end{proof}

\begin{theorem} \label{thm:realizestwt}
  The transformations $x_i \tensor h_i$ realize the graph $\GA$ for $n
  \geq 2$ if we take these values for $x$ and $h$:
  \begin{align*}
    x_1 &= e_1 - e_2 & h_1 &= f_1 + f_2, \\
    x_i &= e_{i - 1} + e_i & h_i &= f_{i - 1} - f_i & \text{for $1
      < i \leq n$.}
  \end{align*}
  The Lie algebra $\LL = \langle x_i \tensor h_i \ranLie$ is $\sll_n$.
\end{theorem}
\begin{proof}
  It is easy to check that the transformations $x_i \tensor h_i$
  realize the graph.  In order to prove that they generate $\sll_n$,
  consider the group $G = \langle 1 + t (x_i \tensor h_i) \mid t \in
  \FF, 1 \leq i \leq n \ranGp$, of which $\LL$ is the Lie algebra. It
  is clear that the orbit of $G$ on $x_1$ spans $V$. Hence, using
  Theorem~\ref{thm:camhall}, either $G = \Sp (V, B)$ or $G = \SL (V)$.
  Hence $\LL$ is either $\spp_n$ or $\sll_n$. Now consider $\exp (2
  \ad x_3 \tensor h_3) (x_1 \tensor h_1) = (x_1 - 2 x_3) \tensor (h_1
  - 2 h_3)$. It forms an extraspecial pair with $x_2 \tensor h_2$,
  because $(h_1 - 2 h_3) (x_2) = 0 \not= h_2 (x_1 - 2 x_3)$. By
  Corollary~\ref{cor:sp-no-exspp}, $\LL = \sll_n$.
\end{proof}

\subsection{Siegel transvections}\label{ss:siegel-transvections}
In order to realize the two orthogonal types of algebras we use Siegel
transvections. As an equivalent of Theorem~\ref{thm:camhall} we will use the
main theorem of Steinbach's paper~\cite{cit:stein97}, dealing with Siegel
transvection groups in a similar way to how Theorem~\ref{thm:camhall} deals
with transvection groups. Steinbach's main theorem is reprinted here in a
weaker form as Theorem~\ref{thm:anja}. In this form, Theorem~\ref{thm:anja}
also follows from the work by Liebeck and Seitz \cite{cit:LS94}.

Let the dimension of $V$ be $2n$ or $2n - 1$. Let $B$ be a
nondegenerate orthogonal bilinear form on $V$; we will denote the
corresponding matrix by $B$ as well. We may assume that $B$ is either
$\bigl(\begin{smallmatrix} 0 & I \\ I & 0
\end{smallmatrix} \bigr)$ or $\Bigl(\begin{smallmatrix}0 & I & 0 \\
  I & 0 & 0 \\ 0 & 0 & 2\end{smallmatrix} \Bigr)$.  Let $u, v \in V$
be linearly independent with $B (u, u) = B (u, v) = B (v, v) = 0$; in
other words, $u$ and $v$ span an \emph{isotropic line}\index{isotropic
  line}. When we speak of lines in an orthogonal
space we mean isotropic lines. Then
\[
S_{u,v}\colon V \to V,\quad x \mapsto x + B (u, x) v - B (v, x) u
\]
is known as the \emph{Siegel transvection}\index{Siegel transvection} determined by $u$ and
$v$.  If $S_{u,v}$ is a Siegel transvection, then we will call the map
$T_{u,v} = S_{u,v} - 1$ an \emph{infinitesimal Siegel transvection}\index{infinitesimal Siegel transvection}.
Note that $T_{u,v}$ is determined up to scalar multiples by the
projective line containing $u$ and $v$. We call the group $\langle 1 +
t T_{u,v} \mid t \in \FF \ranGp$ a \emph{Siegel transvection
  group}\index{Siegel transvection group}.
\begin{lemma}
  $\oo (V)$ is spanned by infinitesimal Siegel transvections.
\end{lemma}
\begin{lemma}
  Infinitesimal Siegel transvections are extremal elements of $\oo
  (V)$.
\end{lemma}
\begin{lemma}
For all $t \in \FF$ and $u, v, w, x \in V$, we have
$\exp (t \ad T_{u,v}) T_{w, x} =  T_{w + t T_{u,v}w, x + t T_{u,v} x}$.
\end{lemma}
In the following lemma, we will let ``point'' refer to projective
points and ``line'' to projective isotropic lines.
\begin{lemma}
  Let $\LL = \langle T_{u, v}, T_{w, x} \ranLie$, where $T_{u,v}$ and $T_{w, x}$ are two infinitesimal Siegel
  transvections.
  Then the isomorphism class of $\LL$ depends on the geometrical
  configuration of the lines $\ell = \langle u, v \ranF$ and $m =
  \langle w, x\ranF$, as follows:
  \begin{itemize}
  \item If $\ell = m$, then $\LL$ is one-dimensional.
  \item If $\ell \cap m$ is a point, then $\LL$ is two-dimensional.
  \item If each point on $\ell$ is collinear to each point on $m$,
    then $\LL$ is two-dimensional.
  \item If exactly one point on $\ell$ is collinear with all of $m$
    and exactly one point on $m$ is collinear with all of $\ell$, then
    $\LL$ is isomorphic to the Heisenberg algebra.
  \item If every point on $\ell$ is collinear with exactly one point
    of $m$ and every point on $m$ is collinear with exactly one point
    of $\ell$, then $\LL$ is isomorphic to $\sll_2$.
  \end{itemize}
  These are all cases.
\end{lemma}
In the case where $\LL$ is isomorphic to the Heisenberg algebra, we
say that $\ell$ and $m$ form an \emph{extraspecial
  pair}\index{extraspecial pair}, just as in the previous section.

We focus our attention on the Siegel transvection groups now. Since
the theorem of Steinbach is somewhat more involved than
Theorem~\ref{thm:camhall}, we need to introduce some notation and
terminology first.

Let $Y = \Omega (V, B)$ be the commutator subgroup of the orthogonal
group $O (V, B)$. Then by~\cite[section~1.1]{cit:stein97} $Y$ is
generated by all Siegel transvection subgroups of $O (V, B)$. Let $G
\not= 1$ be a subgroup of $Y$ generated by some of the Siegel
transvection groups. For each Siegel transvection group $A \subseteq
Y$, let $A^0 = A \cap G$; then either $A^0 = 1$ or $A^0 = A$ -- this
can easily be seen in the Lie algebra.  This situation is a special
case of the situation in the main theorem of~\cite{cit:stein97},
reprinted here in a weaker version as Theorem~\ref{thm:anja}, which
tells us to which isomorphism classes such a group $G$ can belong.
Following~\cite{cit:stein97}, we use the following abbreviations for
different situations:
\begin{description}
\item[(O)] $G$ is an orthogonal group over the same vector space,
  viewed as a vector space over a smaller field. Also, if $B$ has
  maximal Witt index over a vector space of even dimension, then it
  contains the orthogonal group over a vector space of one dimension
  less where $B$ has a one-dimensional radical.
\item[(E)] $G$ is a special case related to triality and $\dim V = 8$.
\item[(ND)] $G$ is the special linear group on a vector space $V'$
  such that $V$ has double the dimension of $V'$ and $B$ has maximal
  Witt index. $V$ is the direct sum of the space $V'$ where $G$ acts
  naturally and a space where it acts dually.
\item[(I)] $G$ is a unitary group.  The unitary space is then regarded
  as orthogonal space over the fixed field of the involutory
  automorphism.
\item[(IQ)] $G$ is a unitary group over a quaternion division ring
  instead of a field.
\item[(KC)] $G$ is a special case related to the Klein correspondence
  and $\dim V \in \{5, 6\}$.
\end{description}
\begin{theorem} \label{thm:anja}
  Let $G\not= 1$ be a subgroup of $Y = \Omega (V, q)$ generated by
  Siegel transvection groups $S_{u_i, v_i}$. Suppose these conditions
  are satisfied:
  \begin{enumerate}
  \item[(H1')] If $A, B$ are Siegel transvection subgroups of $Y$
    that intersect $G$ nontrivially, and $\langle A, B \ranGp =
    \textrm{SL}_2 (\FF)$, then $\langle A^0, B^0 \ranGp = \textrm{SL}_2
    (\FF)$.
  \item[(H2)] All nilpotent normal subgroups of $G$ are contained in
    $Z(G)$ and in the commutator group $G'$, and it is not possible to
    decompose the set $\Sigma$ of Siegel transvection subgroups into two
    nonempty parts $\Sigma_1, \Sigma_2$ such that all groups from
    $\Sigma_1$ commute with all groups from $\Sigma_2$.
  \item[(H3)] There are three pairwise distinct commuting Siegel
    transvection subgroups of $G$ such that there exist Siegel
    transvection subgroups $T$ of $G$ commuting with exactly two of
    these three subgroups.
  \end{enumerate}
  Then one of the situations~(O), (E), (ND), (I), (IQ) and~(KC) holds.
\end{theorem}
\begin{theorem} \label{thm:group-is-orthogonal}
  Suppose all conditions from Theorem~\ref{thm:anja} hold and the
  following extra conditions hold as well.
  \begin{enumerate}[(X1)]
  \item The dimension of $V$ is greater than $8$.
    \label{enum:sufficientlybig}
  \item There are extraspecial pairs in $G$. \label{enum:extraspecial}
  \item There is no $G$-invariant decomposition of $V$ into two
    equal-dimensional subspaces, such that each of these subspaces
    intersects every line corresponding to a Siegel transvection
    subgroup nontrivially. \label{enum:notdirectsum}
  \item One of the Siegel transvection subgroups is isomorphic to
    $\FF^+$. \label{enum:fullfield}
  \item $\langle u_i, v_i \ranF = V$. \label{enum:fullspace}
  \end{enumerate}
  Then $G = Y$.
\end{theorem}
\begin{proof}
  Extra condition~(X\ref{enum:sufficientlybig}) shows that we are not in
  situations~(E) or~(KC).  By extra condition~(X\ref{enum:extraspecial}),
  we must be in situation~(ND) or~(O). By extra
  condition~(X\ref{enum:notdirectsum}), we are not in situation
  (ND). Hence, we are in situation~(O): $G$ is an orthogonal group.
  Then condition~(X\ref{enum:fullfield}) shows that the field is all of
  $\FF$. Finally, because of extra condition~(X\ref{enum:fullspace}), $G$
  is all of $Y$.
\end{proof}

We define the basis of $V$ in an order corresponding to the matrix of
$B$, as follows.  Let $k = n$ for $\Dd$ and $k = n - 1$ for $\Bb$. The
basis of $V$ consists first of vectors $\{e_i\}_{i = 1}^k$ spanning a
maximal isotropic subspace, then of vectors $\{f_i\}_{i = 1}^k$
spanning a maximal isotropic complement to $\langle e_i \ranF$ and
with $B (e_i, f_j) = \delta_{i,j}$, and if $n$ is odd, finally a
vector $g$ with $B (g, g) = 2$.  We can interpret the vectors and
linear functionals from Section~\ref{ss:transvections}, and in
particular from Theorem~\ref{thm:realizestwt}, in this context as
well: $\{f_i\}$ is still the dual basis of $\{e_i\}$. So we
find the well-known isomorphic copy of $\sll (W)$ in $\oo (V)$, where $W$ is
a maximal isotropic subspace; the isomorphism is
determined by sending $x \tensor h \in \sll (W)$ to $T_{x,h} \in \oo
(V)$.  We use the generating elements of $\sll (W)$ in finding those
for $\oo (V)$, but with extra parameters for which the necessity will
become apparent in Section~\ref{s:alwaysours}.
\begin{theorem} \label{thm:realizestwotr}
  Let $\alpha, \beta \in \FF$ and write $\kappa$ for $\sqrt{1 +
    \beta}$. Let
  \[
  \lambda =
  \begin{cases}
    \frac{\alpha}{\alpha + 2}, & \text{if $n$ is odd,} \\
    -\frac{\alpha \kappa}{(\alpha + 2)(1 + \beta + \kappa)}, &
    \text{if $n$ is even.}
  \end{cases}
  \]
  Suppose that $(\alpha + 2)\beta(\beta + 1)\not= 0$ and that $\lambda
  (2 - \beta + \lambda \beta) \ne 1$. The transformations
  $T_{u_i,v_i}$ for $i \leq n$ realize the graph $\GD$ for $n \geq 5$
  if we take these values for $u$ and $v$:
  \begin{align*}
    u_1 &= e_1 - e_2 & v_1 &= f_1 + f_2 + \alpha \tilde f, \\
    u_i &= e_{i - 1} + e_i & v_i &= f_{i - 1} - f_i & \text{for $1
      < i < n$,} \\
    u_n &= e_{n - 2} + \beta f_{n - 1} + e_n & v_n &= f_{n - 2} + e_{n - 1} -
    (1 + \beta) f_n,
  \end{align*}
  where
  \begin{align*}
    \tilde f ={}& (0, 0, 1, 1, \dotsc, 1, 0 \mid 0, 0, 1,
    -1, \dotsc, -1, -1 - \beta), && \text{if $n$ is odd,} \\
    \tilde f ={}& \frac{1}{1 + \beta + \kappa}(0, 0, -\kappa, -\kappa,
    \dotsc, -\kappa, 1 \mid \\ & 0, 0, 1 + \beta + \kappa, -1 - \beta -
    \kappa, \dotsc, 1 + \beta + \kappa, (\beta + 1)(\kappa + 1)), &&
    \text{if $n$ is even;}
  \end{align*}
  where we first write the coefficients of $e_i$, then a bar
  ($\mid$), then those of $f_i$ and finally, if $n$ is odd, the coefficient of $g$.  Then $\LL = \langle T_{u_i,v_i}
  \ranLie$ is $\oo_{2n}$.
\end{theorem}
Note that $B (\tilde f, \tilde f) = 0$, and
\begin{align*}
  B (\tilde f, u_i) &=
  \begin{cases}
    0, & \text{if $i \not= 3$,} \\
    1, & \text{if $i = 3$;} \\
  \end{cases} \\
  B (\tilde f, v_i) &=
  \begin{cases}
    0, & \text{if $i \not= 3$,} \\
    -1, & \text{if $i = 3$ and $n$ is odd,} \\
    \frac{\kappa}{1 + \beta + \kappa}, & \text{if $i = 3$ and $n$ is
      even.}
  \end{cases}
\end{align*}
Hence $\lambda = - \frac{\alpha}{\alpha + 2} B (\tilde f, v_3)$.

We will prove this theorem using
Theorem~\ref{thm:group-is-orthogonal}.  We will first state a similar
theorem for the Lie algebra of type $\Bb$,
Theorem~\ref{thm:realizestrspt}, then state and prove some additional
lemmas necessary to show that the conditions of
Theorem~\ref{thm:group-is-orthogonal} hold, and finally prove
Theorems~\ref{thm:realizestwotr} and~\ref{thm:realizestrspt} on
page~\pageref{proof:twotrspt}.

\begin{theorem} \label{thm:realizestrspt}
  Let $\gamma \in \FF$ be such that $\gamma (\gamma + 1) \not= 0$. The
  transformations $T_{u_i,v_i}$ for $i \in I = \{1, \dotsc, n\}$
  realize the graph $\GB$ for $n \geq 5$ if we
  take these values for $u$ and $v$:
  \begin{align*}
    u_1 &= e_1 - e_2 & v_1 &= f_1 + f_2, \\
    u_i &= e_{i - 1} + e_i & v_i &= f_{i - 1} - f_i & \text{for $1
      < i < n$,} \\
    u_n &= \gamma e_{n - 2} + f_{n - 2} + \gamma e_{n - 1} - f_{n - 1} &
    v_n &= e_{n-2} - f_{n-2} + (1 - \gamma) e_{n - 1} + g,
  \end{align*}
  Then $\LL = \langle T_{u_i, v_i} \ranLie$ is $\oo_{2n - 1}$.
\end{theorem}

The proof will be similar for $\Dd$ and $\Bb$, so let $\LL$ be one of
the two algebras defined in Theorems~\ref{thm:realizestwotr}
and~\ref{thm:realizestrspt} and let $u_i$ and $v_i$ be the
corresponding vectors. Let $G$ be the group generated by the Siegel
transvection groups $S_{u_i, v_i}$. We denote by $\Sigma$ the orbit of
$G$ on the lines $\langle u_i, v_i \ranF$, or alternatively, on the
projective infinitesimal Siegel transvections $\FF T_{u_i, v_i}$, or
alternatively, on the Siegel transvection subgroups $\langle 1 + \mu
T_{u_i, v_i} \mid \mu \in \FF \ranGp$; bijections between these 
three sets are given in Lemma \ref{lm:obvious-bijections}.
Before we start checking the conditions of
Theorem~\ref{thm:group-is-orthogonal}, we will first verify that the
action of $G$ on these three classes of objects is equivalent and that
$G$ is transitive on $\Sigma$. This is asserted by the following two
lemmas, which are readily proven.
\begin{lemma}\label{lm:obvious-bijections}
  The action of $1 + \mu T_{w,x}$ on $\oo (V)$ by conjugation from
  the left is the same as the natural action of $\exp (\mu \ad
  T_{w, x})$. Furthermore, the bijections sending $\FF T_{u, v} \in
  \Sigma$ to $\langle u, v \ranF$ and $\langle 1 + \mu T_{u, v}
  \mid \mu \in \FF \ranGp$ commute with the action of $G$.
\end{lemma}
\begin{lemma}
  All $T_{u_i, v_i}$ are in one orbit under $G$.
\end{lemma}
\begin{lemma} \label{lm:extraspecial-pair}
  $\Sigma$ contains an extraspecial pair.
\end{lemma}
\begin{proof}
  We take $\langle u_2, v_2 \ranF$ as the first line.  For $\Dd$, the
  second line is (the line corresponding to)
  \[
  \exp (2 \ad T_{u_3, v_3}) T_{u_1, v_1} = T_{u_1 + 2 u_3, v_1 + 2 (1 +
    \alpha) v_3 + 2 (\alpha + 2) \lambda u_3};
  \]
  for $\Bb$, we specialize from these values by setting $\alpha =
  \lambda = 0$ and obtain as second line (the line corresponding to)
  \[
  \exp (2 \ad T_{u_3, v_3}) T_{u_1, v_1} = T_{u_1 + 2 u_3, v_1 + 2 v_3}.
  \]
\end{proof}
\begin{lemma} \label{lm:sigma-contains-sln}
  $\Sigma$ contains a copy of the set of infinitesimal transvections
  in $\sll_{n - 1}$.
\end{lemma}
\begin{proof}
  We will exhibit two sets of generators of different isomorphic
  copies of $\sll_{n-1}$. Let $U$ be the vector space spanned by $e_i$
  for $i < n$. If we are studying $\Dd$, then define
  \begin{equation} \label{eq:h_idef} v_i' =
    \begin{cases}
      v_i - \lambda u_i, & i \in \{1, n\} \\
      v_i + \lambda u_i, & 1 < i < n;
    \end{cases}
    \qquad\qquad h_i\colon U \to \FF, x \mapsto - B (v_i, x);
  \end{equation}
  for $\Bb$, substitute $\lambda = 0$ (whence $v_i' = v_i$), but let
  $v_n' = u_n - \gamma v_n$. Additionally, for both $\Bb$ and $\Dd$
  define $u_i' = u_i$, except that $u_n' = \frac{1}{1 + \gamma}(u_n +
  v_n)$ for $\Bb$.  Then $B (u_i', u_j') = B (v_i', v_j') = 0$ for all
  $\{i, j\} \not= \{n-1, n\}$.  Furthermore, $u_i \tensor h_i =
  T_{u_i, v_i} = T_{u_i', v_i'}$ as linear transformations, and $U^*
  \coloneqq \langle h_i \mid i < n \ranF$ is a vector space dual to
  $U$ where the duality is provided by the form $B$. The infinitesimal
  transvections $u_i \tensor h_i$ generate $\sll (U)$ because of the
  same arguments that prove Theorem~\ref{thm:realizestwt}: the centres
  and axes of the transvections span $U$ and its dual, respectively,
  and there is an extraspecial pair (both lines from
  Lemma~\ref{lm:extraspecial-pair} are in $\sll (U)$). The group
  generated by all of the corresponding transvections is a subgroup of
  $G$, and it is transitive on all transvections in $\sll (U)$. So in
  particular, the group of Siegel transvections is transitive on all
  infinitesimal Siegel transvections in $\sll (U)$.

  Similarly, we can define $\tilde U = \langle u_1', u_2', \dotsc,
  u_{n - 2}', u_n' \ranF$ and ${\tilde U}^* = \langle v_1', v_2',
  \dotsc, v_{n-2}', v_n' \ranF$, on which the transvections generate
  $\sll (\tilde U)$.
\end{proof}
For $\Bb$, there is a nontrivial
intersection between $U + \tilde U$ and $U^* + \tilde U^*$: since
\[
\gamma (1 + \gamma) u_n' + v_n' = (1 + \gamma) u_n = (1 +
\gamma)(\gamma u_{n - 1} + v_{n - 1}) = \gamma (1 + \gamma) u_{n -
  1}' + (1 + \gamma) v_{n - 1}',
\]
the intersection is spanned by $u_n' - u_{n - 1}'$, which is $
\gamma (1 + \gamma)$ times $(1 + \gamma) v_{n - 1}' - v_n'$.

By the previous lemma, $\Sigma$ contains an isomorphic copy of
$\sll_4$, which proves the following lemma:
\begin{lemma} \label{lm:ta-tb-tc-td}
  $\Sigma$ contains a $4$-tuple $(T_a, T_b, T_c, T_d)$ of projectively distinct infinitesimal
  Siegel transvections  such that $T_c$ and
  $T_d$ do not commute, but every other pair does.
\end{lemma}
\begin{lemma} \label{lm:double-centralizer}
  Let $T_a = T_{u_a, v_a} \in \Sigma$ and $T_b =
  T_{u_b, v_b} \in \Sigma$ satisfy $C_\Sigma (T_a) = C_\Sigma (T_b)$.
  Then $T_a = T_b$.
\end{lemma}
\begin{proof}
  Since $G$ is transitive on $\Sigma$, we may assume that $T_a =
  T_{u_2, v_2} = T_{u_2', v_2'}$. Let us denote the subspace of $V$
  perpendicular with respect to the bilinear form $B$ to a vector $u
  \in V$ by $u^\perp$, and similarly, let us denote the subspace of
  $V$ perpendicular to a subspace $S$ of $V$ by $S^\perp$.

  Recall from Lemma~\ref{lm:sigma-contains-sln} the definitions of $U$
  and $U^*$. Pick a nonzero vector $u \in U$ which is perpendicular to
  $v_2'$, and let nonzero $v \in U^*$ be perpendicular to $u$ and to
  $u_2'$.  Then the infinitesimal transvection $u \tensor v = T_{u,
    v}$ is in $\sll (U)$, hence in $\Sigma$, and it commutes with
  $T_{u_2', v_2'}$. So $T_b$ should also commute with it. Then
  $\langle u_b, v_b \ranF$ either intersects $\langle u, v \ranF$, or
  $u_b$ and $v_b$ are both perpendicular to $u$ and $v$.  So if
  $\langle u_b, v_b \ranF$ does not intersect all $\langle u, v \ranF$
  for fixed $u$ and all nonzero $v \in S \coloneqq u^\perp \cap
  u_2'^\perp \cap U^*$, then $u_b$ and $v_b$ are perpendicular to $u$.
  Let us consider the case where $\langle u_b, v_b \ranF$ intersects
  every such $\langle u, v \ranF$. We will show that $u_b$ and $u_v$
  are perpendicular to $u$ in this case as well. $S$ has codimension
  $1$ or $2$ in $U^*$ of dimension $n - 1$, so its dimension is at
  least $2$. If $\dim S > 2$, then $\langle u_b, v_b \ranF$ must
  contain $u$ to intersect every $\langle u, v \ranF$. In that case,
  $u_b$ and $v_b$ are certainly perpendicular to $u$. Hence assume
  $\dim S = 2$ and $u \not\in \langle u_b, v_b \ranF$. Then for
  different lines $\langle u, v \ranF$, the intersection with $\langle
  u_b, v_b \ranF$ is different. So the intersections span all of
  $\langle u_b, v_b \ranF$. In particular, $u_b$ and $v_b$ themselves
  are on lines $\langle u, v \ranF$. Thus both are perpendicular to
  $u$.

  We see that $u_b$ and $v_b$ are perpendicular to all of $v_2'^\perp
  \cap U$. Similarly, they are perpendicular to $v_2'^\perp \cap
  \tilde U$; that is, they are perpendicular to $S_v \coloneqq
  v_2'^\perp \cap (U + \tilde U)$. Dually, we see that $u_b$ and $v_b$
  are perpendicular to $S_u \coloneqq u_2'^\perp \cap (U^* + \tilde
  U^*)$.

  For $\Bb$, the intersection of $U + \tilde U$ and $U^* + \tilde U^*$
  is spanned by $u_n' - u_{n-1}'$, which is perpendicular to both
  $u_2'$ and $v_2'$; for $\Dd$, the intersection of $U + \tilde U$
  with $U^* + \tilde U^*$ is empty. So $S_u + S_v$ is a $2n -
  2$-dimensional space for $\Dd$ and to a $2n - 3$-dimensional space
  for $\Bb$. Since the form is nondegenerate, there is in both cases
  only a $2$-dimensional space of vectors perpendicular to $S_u +
  S_v$.  This space is $\langle u_2', v_2'\ranF$. Hence $u_b$ and
  $v_b$ are in $\langle u_2', v_2' \ranF$.
\end{proof}
\begin{lemma} \label{lm:graph-f-sigma}
  The graph $F (\Sigma)$ with vertex set $\Sigma$ and where two
  infinitesimal Siegel transvections are adjacent if they generate an
  algebra isomorphic to $\sll_2$, is connected.
\end{lemma}
\begin{proof}
  According to Lemma~2.13 of~\cite{cit:timme01}, if $\card \Sigma >
  1$, then $F (\Sigma)$ is connected if and only if $\Sigma$ is a
  conjugacy class in $G$ and $F (\Sigma)$ has an edge. Both of these
  conditions are fulfilled.
\end{proof}
\begin{lemma} \label{lm:quasisimple-group}
  $G$ is a quasisimple group.
\end{lemma}
\begin{proof}
  We use Lemma~2.14 of~\cite{cit:timme01}. A weaker version of it
  states that if the following conditions are satisfied:
  \begin{itemize}
  \item $\card \Sigma > 1$;
  \item the graph $F (\Sigma)$ from Lemma~\ref{lm:graph-f-sigma} is
    connected;
  \item there exists no pair $A \not= C \in \Sigma$ with $C_\Sigma (A)
    = C_\Sigma (B)$;
  \item $\Sigma$ contains an extraspecial pair;
  \item the elements of $\Sigma$ correspond to extremal Lie algebra
    elements;
  \end{itemize}
  then $G$ is a quasisimple group.  Lemmas~\ref{lm:extraspecial-pair},
  \ref{lm:double-centralizer} and~\ref{lm:graph-f-sigma} show that the
  three nontrivial conditions are fulfilled.
\end{proof}

The technical proof of the following lemma uses the previous lemmas and 
has been omitted here for brevity. It can be found in \cite{cit:postm07}.

\begin{lemma} \label{lm:condition-h2-x3}
  $G$ satisfies conditions (H2) from Theorem~\ref{thm:anja} and
 	(X3) from Theorem~\ref{thm:group-is-orthogonal}.
\end{lemma}

\begin{proof}[Proof of Theorems~\ref{thm:realizestwotr}
  and~\ref{thm:realizestrspt}] \label{proof:twotrspt}
  We intend to apply Theorem~\ref{thm:group-is-orthogonal}, so we will
  need to show that its conditions hold.
  \begin{itemize}
  \item Condition~(H1') follows from the fact that for every Siegel
    transvection subgroup $A$, we have $A^0 = 1$ or $A^0 = A$.
  \item Condition~(H2) follows from Lemma~\ref{lm:condition-h2-x3}.
  \item Condition~(H3) follows from Lemma~\ref{lm:ta-tb-tc-td}.
  \item Condition~(X\ref{enum:sufficientlybig}) is clearly satisfied.
  \item Condition~(X\ref{enum:extraspecial}) follows from
    Lemma~\ref{lm:extraspecial-pair}.
  \item Condition~(X\ref{enum:notdirectsum}) follows from Lemma~\ref{lm:condition-h2-x3}.
  \item Condition~(X\ref{enum:fullfield}) is true for every Siegel
    transvection subgroup.
  \item Condition~(X\ref{enum:fullspace}) is clearly satisfied. \qedhere
  \end{itemize}
\end{proof}

\section{The algebras nearly always correspond to these realizations}\label{s:alwaysours}

In this section, we show that a Lie algebra $\LL$ realizing one of the
graphs~$\GA$, $\GB$, $\GC$ and~$\GD$, is in the generic case a
quotient of the realization $\MM$ we found in the previous section. In
order to see this, we find different sets of generators. For types
$\Bb$ and $\Dd$, we will need the parameters~$\alpha$, $\beta$
and~$\gamma$ from Theorems~\ref{thm:realizestwotr}
and~\ref{thm:realizestrspt} to have sufficient degrees of freedom to
be able to make the sets of generators of $\LL$ and $\MM$ match
up.

Since $\MM$ is simple in most cases, it will follow that $\LL$ and
$\MM$ are isomorphic. The only exception is $\Aa$ if $p \mid n$, as is
well known from literature (see e.g.~\cite{cit:humph72}
and~\cite{cit:strad04}).

\begin{theorem} \label{thm:linearalways}
  Let $\LL = \langle x_1, \dotsc, x_n \ranLie$ be any Lie algebra
  realizing the graph~$\GC$ from Figure~\ref{fig:linear}. Suppose that
  $n$ is even and that $f (x_i, x_{i + 1}) \ne 0$ for all $i$. Then
  $\LL$ is isomorphic to $\spp_n$.
\end{theorem}
\begin{proof}
  Let $\MM$ be the realization from Theorem~\ref{thm:realizeslinear}.
  By that Theorem, $\MM = \spp_n$.  Denote the generators of $\MM$
  realizing $\GC$ by $z_i$. Scale $x_i$ such that the extremal form is
  identical on both sets of generators. Then the map $\sigma\colon \MM
  \to \LL$ mapping each of the monomials $y_{k, m}$ in $z_i$ to the
  same monomial in $x_i$, is a Lie algebra homomorphism by
  Lemma~\ref{lm:f-C}. Hence $\LL$ is a quotient of $\MM$. But since
  $\MM$ is simple, $\LL \cong \MM$.
\end{proof}

Already for Lie algebras of type~$\Aa$, we need more degrees of
freedom than can be obtained by just scaling the generators. The
following lemma, which is based on Section~5 of~\cite{cit:cosuw01},
will be sufficient.
\begin{lemma} \label{lm:fixtriangle}
  Let $\pi, \rho, \sigma \in \FF$ all be nonzero. Let $x$, $y$, $z$ be
  extremal elements of $\LL$ such that $f (x, yz)^2 \not= 2 f (x, y) f
  (x, z) f (y, z)$ and $f (x, y) \not= 0 \not= f (y, z)$, and such
  that $x$ and $y$ commute with a set $S$ of elements of $\LL$. We can
  find extremal elements $x'$, $y'$ and $z'$ with the following
  properties:
  \begin{itemize}
  \item $\langle x, y, z \ranLie = \langle x', y', z' \ranLie$,
  \item $x'$ and $y'$ commute with $S$,
  \item $(f (x', y'), f (x', z'), f (y', z'), f (x', y'z')) = (\pi, \rho,
    \sigma, 0).$
  \end{itemize}
\end{lemma}
\begin{proof}
  Let $s = f (x, yz) / (f (x, y) f (y, z))$.  Let $\hat{x} = \exp (s
  \ad y) (x)$. Then $\langle x, y, z \ranLie = \langle \hat x, y, z
  \ranLie$ (since $\exp (- s \ad y) (\hat x) = x$) and
  \begin{align*}
    f (\hat x, y) &{}= f (x, y), \\
    f (\hat x, z) &{}= f (x, z) - \frac{f (x, yz)^2}{2 f (x, y) f (y, z)}, \\
    f (\hat x, yz) &{}= 0.
  \end{align*}
  Note that $f (\hat{x}, z) \not= 0$. We drop the circumflex from now
  on and use $x$ to denote $\hat x$.  Scale $x$, $y$ and $z$ to obtain
  $\tilde{x}$, $\tilde{y}$ and $\tilde{z}$:
  \begin{align*}
    \tilde{x} &{}= \sqrt{\frac{\pi \rho f (y, z)}{\sigma f (x, y) f (x, z)}} x, \\
    \tilde{y} &{}= \sqrt{\frac{\pi \sigma f (x, z)}{\rho f (x, y) f (y, z)}} y, \\
    \tilde{z} &{}= \sqrt{\frac{\rho \sigma f (x, y)}{\pi f (x, z) f (y, z)}} z.
  \end{align*}
  Now all values of $f$ are as intended, possibly up to a factor of
  $-1$, and $f (\tilde x, \tilde y) f (\tilde x, \tilde z) f (\tilde
  y, \tilde z) = \pi \rho \sigma$. So either all values of $f$ are
  exactly as intended (including sign), in which case we are done, or
  exactly two of them need their sign changed; say $f (\tilde x,
  \tilde y)$ and $f (\tilde x, \tilde z)$.  Then let $x' = - \tilde
  x$, $y' = \tilde y$ and $z' = \tilde z$.

  $y'$  commutes with $S$, since it is merely a scaled version of
  $y$. By the Jacobi identity, $[x, y]$ commutes with $S$, as well;
  hence $x' \in \langle x, [x, y], y \ranF$ commutes with $S$.
\end{proof}

\begin{theorem} \label{thm:twtalways}
  Let $\LL = \langle x_1, \dotsc, x_n \ranLie$ be a realization of the
  graph $\GA$ in Figure~\ref{fig:twt}. Suppose that the following Zariski-open
  conditions on the extremal form hold:
  \begin{itemize}
  \item $f (x_1, x_2 x_3)^2 \ne 2 f (x_1, x_2) f (x_1, x_3) f
    (x_2, x_3)$,
  \item $f (x_i, x_{i+1}) \ne 0$ for all $i$.
  \end{itemize}
  Then:
  \begin{itemize}
  \item if $\characteristic \FF = p > 0$ and $p \mid n$, then $\LL$ is
    isomorphic to either $\sll_n$ or its simple subalgebra of
    codimension $1$;
  \item otherwise, $\LL$ is isomorphic to $\sll_n$.
  \end{itemize}
\end{theorem}
\begin{proof}
  Let $\MM$ be the realization from Theorem~\ref{thm:realizestwt}. By
  that theorem, $\MM = \sll_n$.  Denote the generators of $\MM$
  realizing $\GA$ by $y_i$. We will exhibit a Lie algebra homomorphism
  from $\MM$ to $\LL$, showing that $\LL$ is a quotient of $\MM$.
  Since $\LL$ cannot be one-dimensional, we then have the desired
  result.

  Apply Lemma~\ref{lm:fixtriangle} with $(\pi, \rho, \sigma) = (1, 1,
  1)$ and $(x, y, z) = (x_1, x_2, x_3)$. We obtain a new set of
  generators $x_i'$ of $\LL$ that still realize $\GA$. Also apply
  Lemma~\ref{lm:fixtriangle} with $(x, y, z) = (y_1, y_2, y_3)$ and
  with the same values for $\pi$, $\rho$ and $\sigma$, obtaining new
  generators $y_i'$. Now for any pair of monomials in $y_1'$, $y_2'$
  and $y_3'$, the extremal form on that pair is equal to the extremal
  form on the corresponding pair of monomials in $x_1'$, $x_2'$ and
  $x_3'$. For $i > 3$, inductively define $x_i' = f (y_{i - 1}, y_i) f
  (x_{i - 1}', x_i)^{-1} x_i$ and $y_i' = y_i$. Now the extremal form
  is equal on all pairs of monomials given by $\psiA (\ff)$, so the
  extremal form is identical. Then the desired Lie algebra
  homomorphism can be obtained as mapping $x_i'$ to $y_i'$.
\end{proof}
\begin{theorem} \label{thm:twotralways}
  Let $\LL = \langle x_1, \dotsc, x_n \ranLie$ be a realization of the
  graph~$\GD$ in Figure~\ref{fig:twotr}. Then $\LL$ is isomorphic to
  $\oo_{2n}$ if the values of the extremal form satisfy a number of
  Zariski-open conditions.
\end{theorem}
These open conditions can be found as follows. If $n$ is odd,
\begin{itemize}
\item we have the condition $f (x_1, x_{3 \up n - 2} x_{n \down 2})
  \ne 8$,
\item furthermore, we define $\alpha$ by Eq.~\eqref{eq:alpha-odd} and
  $\beta$ by Eq.~\eqref{eq:beta-odd};
\end{itemize}
if $n$ is even,
\begin{itemize}
\item we define $\alpha$ by Eq.~\eqref{eq:alpha-even} and $\beta$ by
  Eq.~\eqref{eq:beta-even};
\end{itemize}
then the (other) conditions are
\begin{itemize}
\item $f (x_1, x_2 x_3)^2 \ne 2 f (x_1, x_2) f (x_1, x_3) f
  (x_2, x_3)$,
\item $f (x_n, x_{n-1} x_{n-2})^2 \ne 2 f (x_n, x_{n-1}) f (x_n,
  x_{n-2}) f (x_{n-1}, x_{n-2})$,
\item $f (x_i, x_{i+1}) \ne 0$ for all $i$,
\item $(\alpha + 2) \beta (\beta + 1) \ne 0$,
\item $\lambda (2 - \beta + \lambda \beta) \ne 1$.
\end{itemize}
Note that if $n$ is odd and $f (x_1, x_{3 \up n - 2} x_{n \down 2}) =
8$, then $\LL$ is also isomorphic to $\oo_{2n}$, but under slightly
different conditions. This can be seen from the proof below.
\begin{proof}
  Let $\MM$ be the realization we defined in
  Theorem~\ref{thm:realizestwotr}, for some values of $\alpha$ and
  $\beta$ which we will choose later. We will make some changes to the
  sets of generating elements that do not change the algebra generated
  by these elements and then claim that $f$ is identical on the two.
  This shows that $\LL$ is isomorphic to a quotient of $\Dd$. Since
  $\Dd$ is simple, this quotient is $\Dd$ itself.

  Let $x_i$ be the $i$th extremal generator of $\LL$ (so the value of
  $x_i$ will change during the rest of this proof).  First, we perform
  the procedure of Lemma~\ref{lm:fixtriangle} above on $x_1$, $x_2$
  and $x_3$ with $(\pi, \rho, \sigma) = (-8, 1, 2)$. Then
  we do the same on $x_{n-2}$, $x_{n-1}$ and $x_n$, this time with
  $(\pi, \rho, \sigma) = (2, 2, 1)$. Now there are $n - 5$ pairs of
  elements left (on the ``connecting line between the two
  triangles''), on which the value of $f$ has not yet been adjusted,
  and additionally $f (x_1, x_{3 \up n-2} x_{n \down 2})$. We assume
  that $f (x_{i - 1}, x_i) \not= 0$ for $4 \leq i \leq n - 3$ and
  scale $x_4$ up to $x_{n - 3}$ such that $f (x_{i - 1}, x_i) = 2$.
  This leaves $f (x_{n - 3}, x_{n - 2})$ and $f (x_1, x_{3 \up n-2}
  x_{n \down 2})$. The values of $f$ other than $f (x_1, x_{3 \up n-2}
  x_{n \down 2})$ are as given in Figure~\ref{fig:f-values}.

  \begin{figure}
    \centering
    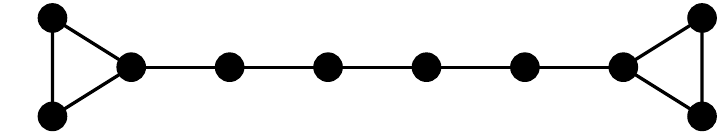
   
    \caption{The values of $f$ on $\LL$. The numbering of the nodes is
      the same as in Figure~\ref{fig:twotr}.}
    \label{fig:f-values}
  \end{figure}

  We now perform a similar procedure on $\MM$. Call the $i$th extremal
  generator $z_i = T_{u_i, v_i}$. The values of $f$ prior to any
  changes are as given in Figure~\ref{fig:f-values-m}. We perform the
  procedure of Lemma~\ref{lm:fixtriangle} on $z_1$, $z_2$ and $z_3$,
  with $(\pi, \rho, \sigma) = (-8, 1, 2)$. Note that $s = \alpha/4$.
  In the first step, $z_1$ becomes
  \[
  z_1 - \frac{\alpha}{4} [z_1, z_2] - \frac{\alpha^2}{4} z_2 =
  T_{u_1 - \frac{\alpha}{2} u_2, v_1 + \frac{\alpha}{2} v_2}.
  \]
  Now
  \begin{align*}
    f(z_1, z_2) &= -8,& f(z_1, z_3) &= \frac{(\alpha + 2)^2}{2},\\
    f (z_2, z_3) &= 2,& f (z_1, z_2 z_3) &= 0,
  \end{align*}
  so $z_1$ and $z_3$ are divided by $(\alpha + 2)/\sqrt{2}$ and
  $z_2$ is multiplied with that same constant.

  \begin{figure}
    \centering 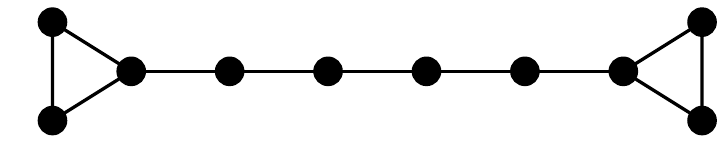
   
    \caption{The values of $f$ on $\MM$, prior to any changes. The
      numbering of the nodes is the same as in
      Figure~\ref{fig:twotr}.}
    \label{fig:f-values-m}
  \end{figure}

  At the other end, we find that applying the procedure of
  Lemma~\ref{lm:fixtriangle} only entails scaling the generators by
  a factor; in particular, for $(\pi, \rho, \sigma) = (2, 2, 1)$, we
  divide $z_{n - 1}$ and $z_n$ by $\sqrt{-2 - 2 \beta}$ and multiply
  $z_{n - 2}$ by the same factor. Finally, to obtain a $2$ for the
  value of $f (z_{i - 1}, z_i)$ where $4 \leq i \leq n - 3$, we
  multiply $z_i$ by $(\alpha + 2)/\sqrt{2}$ for even $i$
  and divide it by that factor for odd $i$. Hence
  \begin{equation} \label{eq:fn}
    f (z_{n - 3}, z_{n - 2}) = (2\alpha + 4)^{(-1)^{n+1}}\sqrt{-1 - \beta}.
  \end{equation}

  Now consider $z = z_{3 \up n - 2} z_{n \down 2}$. The generators
  involved have been scaled, but not changed otherwise. For the
  generators that occur twice, viz.~$z_3, z_4, \dotsc, z_{n - 2}$,
  their scaling factor affects $z$ twice; the scaling factors of the
  other three ($z_2$, $z_{n - 1}$, and $z_n$) have an effect only
  once. In total, $z$ was multiplied by $(\alpha + 2)/\sqrt{2}$ by all
  these scalings, if $n$ is odd, and divided by that constant if $n$
  is even. We will now compute the value of $z$ explicitly. This is
  easier if we temporarily forget all the scaling that occurred; so
  until further notice, we will use the values before scaling of the
  $z_i$.

  With induction it is easy to see that $z_{k \down 2} = (-1)^k
  (T_{u_k, v_2} + T_{u_2, v_k})$ for $3 \leq k < n$. Then
  \begin{multline*}
    z_{n \down 2} = (-1)^{n - 1} ([z_n, T_{u_{n - 1}, v_2}] + [z_n,
    T_{u_2, v_{n-1}}]) \\
    = (-1)^n (T_{u_n, v_2} - \beta T_{v_n, v_2} - T_{u_2, v_n} +
    T_{u_n, u_2} = (-1)^n (T_{u_n, u_2 + v_2} - T_{u_2 - \beta v_2,
      v_n}).
  \end{multline*}
  We can again use induction to see that $z_{k \up n-2} z_{n \down 2} = (-1)^n (T_{u_k, u_2 + v_2} - T_{u_2 - \beta v_2, v_k})$ for $4
  \leq k \leq n - 2$. Finally, we compute
  \begin{multline*}
    z = z_{3 \up n - 2} z_{n \down 2} = (-1)^n([z_3, T_{u_4, u_2 +
      v_2}] - [z_3, T_{u_2 - \beta v_2, v_4}])= \\
    (-1)^n(- T_{u_4, v_3} + T_{u_3, u_2 + v_2} + T_{u_3, u_4} -
    T_{u_2 - \beta v_2, v_3} + T_{u_3, v_4} - \beta T_{v_3, v_4}) = \\
    (-1)^n (T_{u_3, u_2 + v_2 + u_4 + v_4} - T_{u_2 - \beta
      v_2 + u_4 - \beta v_4, v_3}).
  \end{multline*}

  We re-remember the scaling factors and find that
  \begin{multline*}
    f (z_1, z) z_1 = [z_1, [z_1, z]] = \\
     \frac{2(-1)^n}{(\alpha + 2)^2}\left(\frac{\alpha + 2}{\sqrt
        2}\right)^{(-1)^{n + 1}}         [T_{u_1 - \frac{\alpha}{2} u_2, v_1
      + \frac{\alpha}{2} v_2}, [T_{u_1 - \frac{\alpha}{2} u_2, v_1 +
      \frac{\alpha}{2} v_2}, T_{u_3, u_2 + v_2 + u_4 + v_4} - T_{u_2
      - \beta v_2 + u_4 - \beta v_4, v_3}]] \\
    =
    \begin{cases}
      f (T_{u_1 - \frac{\alpha}{2} u_2, v_1 + \frac{\alpha}{2} v_2},
      T_{u_2 - \beta v_2 + u_4 - \beta v_4, v_3} - T_{u_3, u_2 + v_2
        + u_4 + v_4}) z_1, & \text{if $n$ is odd,} \\
      \frac{2}{(\alpha + 2)^2} f([T_{u_1 - \frac{\alpha}{2} u_2, v_1
        + \frac{\alpha}{2} v_2}, T_{u_3, u_2 + v_2 + u_4 + v_4} -
      T_{u_2 - \beta v_2 + u_4 - \beta v_4, v_3}) z_1, & \text{if
        $n$ is even.}
    \end{cases}
  \end{multline*}
  A straightforward computation shows that
  \[
  f (T_{u_1 - \frac{\alpha}{2} u_2, v_1 + \frac{\alpha}{2} v_2},
  T_{u_2 - \beta v_2 + u_4 - \beta v_4, v_3} - T_{u_3, u_2 + v_2 +
    u_4 + v_4}) = 4 \alpha (1 + \beta) + 8.
  \]
  We finish the proof by case distinction.

  \begin{caselist}
  \item Suppose that $n$ is odd. Taking the previous equations
    together with Eq.~\eqref{eq:fn}, we need to choose $\alpha$
    and $\beta$ such that
    \begin{align}
      f (x_1, x_{3 \up n-2} x_{n \down 2}) &{}= 4 \alpha (1 +
      \beta)
      + 8 \label{eq:flong-odd}\\
      \intertext{and} f (x_{n-3}, x_{n-2}) &{}= (2 \alpha +
      4)\sqrt{-1 - \beta}. \label{eq:fshort-odd}
    \end{align}
    If $f (x_1, x_{3 \up n-2} x_{n \down 2}) = 8$, then we choose
    $\alpha = 0$ and $\beta = -1-f (x_{n-3}, x_{n-2})^2/16$,
    otherwise we use Eq.~\eqref{eq:flong-odd} to obtain
    \begin{equation} \label{eq:beta-odd} 1 + \beta = \frac{f (x_1,
        x_{3 \up n-2} x_{n \down 2}) - 8}{4 \alpha};
    \end{equation}
    substituting this into Eq.~\eqref{eq:fshort-odd}, we obtain
    \[
    (\alpha + 2) \sqrt{\frac{8 - f (x_1, x_{3 \up n-2} x_{n \down
          2})}{\alpha}} = f (x_{n-3}, x_{n-2})
    \]
    or
    \[
    \sqrt{\alpha} + \frac{2}{\sqrt{\alpha}} = \frac{f (x_{n-3},
      x_{n-2})}{\sqrt{8 - f (x_1, x_{3 \up n-2} x_{n \down 2})}}.
    \]
    This is a quadratic equation in $\sqrt{\alpha}$ that can easily
    be solved to
    \begin{equation} \label{eq:alpha-odd}
    \alpha = \pm \frac{\left(f (x_{n-3}, x_{n-2}) \pm \sqrt{f
          (x_{n-3}, x_{n-2})^2 + 8 f (x_1, x_{3 \up n-2} x_{n \down
            2}) - 64} \right)^2}{32 - 4 f (x_1, x_{3 \up n-2} x_{n
        \down 2})},
  \end{equation}
    four solutions (some of which may coincide). Then $\beta$ can be
    found from Eq.~\eqref{eq:beta-odd}. Now all values of $f$
    are the same and hence $\LL$ is a quotient of $\MM$. Since $\MM$
    is simple, $\LL = \MM$.
  \item Suppose that $n$ is even. We need to choose $\alpha$ and
    $\beta$ such that
    \begin{align}
      f (x_1, x_{3 \up n-2} x_{n \down 2}) &{}= -\frac{8 \alpha (1 +
        \beta) + 16}{(\alpha + 2)^2} \label{eq:flong-even} \\
      \intertext{and} f (x_{n-3}, x_{n-2}) &{}= (2 \alpha +
      4)\sqrt{-1 - \beta}. \notag
    \end{align}
    The last equation can be written as
    \begin{equation} \label{eq:beta-even}
      \beta = - \frac{f (x_{n-3}, x_{n-2})^2}{2 (\alpha + 2)} - 1.
    \end{equation}
    If we substitute this into Eq.~\eqref{eq:flong-even}, we get
    \begin{equation} \label{eq:alpha-even}
      f (x_1, x_{3 \up n-2} x_{n \down 2}) (\alpha + 2)^3 = 8 \alpha
      f (x_{n-3}, x_{n-2})^2 - 16 \alpha - 32,
    \end{equation}
    which we can solve for $\alpha$ and thus find an explicit value
    for $\beta$. We find an isomorphism again. \qedhere
  \end{caselist}
\end{proof}

\begin{theorem} \label{thm:trsptalways}
  Let $\LL = \langle x_1, \dotsc, x_n \ranLie$ be a realization of the
  graph~$\GB$ in Figure~\ref{fig:trspt}. Then $\LL$ is isomorphic to
  $\oo_{2n - 1}$ if the values of the extremal form satisfy these
  Zariski-open conditions:
  \begin{itemize}
  \item $f (x_1, x_2 x_3)^2 \ne 2 f (x_1, x_2) f (x_1, x_3) f (x_2, x_3)$,
  \item $f (x_i, x_{i + 1}) \ne 0$ for $i \leq n - 2$,
  \item $f (x_{n - 2}, x_n) \ne 0$,
  \item $f (x_1, x_{3 \up n - 2} x_{n \down 2}) \ne 0$,
  \item $f (x_1, x_{3 \up n - 2} x_{n \down 2}) \ne 8\cdot (-1)^{n+1}$.
  \end{itemize}
\end{theorem}
\begin{proof}
  Let $\MM$ be the realization defined in
  Theorem~\ref{thm:realizestrspt}; we will specify the value of
  $\gamma$ later. Let $z_i = T_{u_i, v_i}$ be the $i$th extremal
  generator of $\MM$. In the same manner as in the proof of
  Theorem~\ref{thm:twotralways}, we change $x_i$ such that the values
  of $f (x_i, x_j)$ are equal to those of $f (z_i, z_j)$ for $i, j
  \leq 3$, and such that $f (x_1, x_2 x_3) = f (z_1, z_2 z_3)$. By
  scaling $x_i$ we can assure that $f (x_{i - 1}, x_i) = f (z_{i-1},
  z_i)$ for $i < n$, and by scaling $x_n$ we can assure that $f
  (x_{n-2}, x_n) = f (z_{n-2}, z_n)$.  Then what remains is assuring
  that $f (x_1, x_{3 \up n-2} x_{n \down 2}) = f (z_1, z_{3 \up n-2}
  z_{n \down 2})$.

  Like in the proof of Theorem~\ref{thm:twotralways}, we will do this
  by choosing the value of $\gamma$ appropriately. This requires
  explicitly constructing $z_{3 \up n-2} z_{n \down 2}$. Again like
  before, for $3 \leq k < n$ we find with induction that $z_{k \down 2}
  = (-1)^k (T_{u_k, v_2} + T_{u_2, v_k})$. Then
  \[
  z_{n \down 2} = (-1)^{n-1} ([z_n, T_{u_{n-1}, v_2}] + [z_n, T_{u_2,
    v_{n-1}}]) = (-1)^{n-1} T_{u_n, \gamma u_2 + v_2}.
  \]
  With induction we see that $z_{k \up n-2} z_{n \down 2} = (-1)^{n-1}
  T_{\gamma u_k + v_k, \gamma u_2 + v_2}$ for $3 < k \leq n-2$. Then
  \[
  z_{3 \up n-2} z_{n \down 2} = (-1)^{n-1} [T_{u_3, v_3}, T_{\gamma
    u_4 + v_4, \gamma u_2 + v_2}] = (-1)^{n-1}T_{\gamma u_3 + v_3, \gamma u_2 +
    v_2 + \gamma u_4 + v_4}.
  \]
  Then
  \begin{multline*}
    f (z_1, z_{3 \up n-2} z_{n \down 2}) = 2 (B (u_1, \gamma u_3 +
    v_3) B (v_1, \gamma u_2 + v_2 + \gamma u_4 + v_4) - \\
    B (u_1, \gamma u_2 + v_2 + \gamma u_4 + v_4) B (v_1, \gamma u_3 +
    v_3)) = (-1)^n 8 \gamma.
  \end{multline*}
  So we choose $\gamma$ to be $(-1)^n/8$ times the value of $f (x_1,
  x_{3 \up n-2} x_{n \down 2})$. Then $f$ is identical on $\LL$ and
  $\MM$, so $\LL$ is a quotient of $\MM$; since $\MM$ is simple, they
  are isomorphic.
\end{proof}

For each of the four infinite families of Chevalley type Lie algebras,
we have given a family of graphs such that a generic Lie algebra
generated by a set of extremal elements realizing such a graph, is
isomorphic to the corresponding Lie algebra.

\section*{Acknowledgements}
We would like to thank Arjeh~Cohen, Hans~Cuypers, and Jan~Draisma for 
various fruitful discussions on the topic. We are also thankful to the 
referee for the thorough review and valuable suggestions.

\def\cprime{$'$}

\end{document}